\documentclass[12pt]{article}
\usepackage{amssymb}
\usepackage{amsmath}
\usepackage{pst-all}
\newtheorem{theorem}{Theorem}

\newtheorem{proposition}[theorem]{Proposition}
\newtheorem{lemma}[theorem]{Lemma}
\newtheorem{corollary}[theorem]{Corollary}
\newtheorem{remark}[theorem]{Remark}

\newtheorem{definition}[theorem]{Definition}
\newcommand{\R}{\mathbb{R}}
\newcommand{\Sf}{\mathbb{S}}
\newcommand{\spa}{\mbox{span\,}}
\newcommand{\hess}{\mbox{\em  Hess\,}}

\newcommand{\grad}{\mbox{grad}}
\newcommand{\po}{{\hspace*{-1ex}}{\bf .  }}
\newcommand{\ii}{isometric immersion }

\newcommand{\sff}{second fundamental form }
\newcommand{\an}{\;\;\;\mbox{and}\;\;\;}

\newcommand{\E}{{\cal E}}
\newcommand{\N}{{\cal N}}
\def\Ral{{\cal R}}
\def\Fes{{\cal F}}
\def\P{{\cal P}}

\def\<{\langle}

\def\>{\rangle}
\def\va{\varphi}
\def\o{\omega}

\def\d{\partial}
\def\bea{\begin{eqnarray*} }
\def\eea{\end{eqnarray*} }
\def\be{\begin{equation} }
\def\ee{\end{equation} }

\def\proof{\noindent{\it Proof: }}
\def\qed{\ifhmode\unskip\nobreak\fi\ifmmode\ifinner\else
\hskip5 pt \fi\fi\hbox{\hskip5 pt \vrule width4 pt
height6 pt  depth1.5 pt \hskip 1pt }}
\begin{document}

\title{The vectorial Ribaucour transformation for submanifolds and applications}
\author{M. Dajczer, L. A. Florit  \& R. Tojeiro}
\date{}
\maketitle
\begin{abstract}
In this paper we develop the vectorial Ribaucour transformation
for Euclidean submanifolds.  We prove a general decomposition
theorem showing that under \mbox{appropriate} conditions the
composition of two or more vectorial Ribaucour transformations is
again a vectorial Ribaucour transformation. An immediate
consequence of this result is the classical permutability of
Ribaucour transformations. Our main application is an explicit
local construction of all Euclidean submanifolds with flat normal
bundle. Actually, this is  a particular case of a more general
result. Namely, we obtain a local explicit construction  of all
Euclidean submanifolds carrying a parallel flat normal subbundle,
in particular of all those  that carry a parallel normal vector
field.  Finally, we describe all submanifolds carrying a Dupin
principal curvature normal vector field with integrable conullity,
a concept that has proven to be crucial in the study of
reducibility of Dupin submanifolds.
\end{abstract}

An explicit construction of all submanifolds with flat normal bundle of
the Euclidean sphere carrying a holonomic net of curvature lines, that is,
admitting principal coordinate systems, was given by Ferapontov in \cite{fe}. The author
points out that his construction ``resembles'' the vectorial Ribaucour
transformation for orthogonal systems developed  in
\cite{lm}. The latter provides a convenient framework for understanding
the permutability properties of the classical Ribaucour transformation.

 This paper grew out as an attempt to better understand the connection
 between those two subjects, as a means of unraveling the geometry
 behind Ferapontov's construction. This has led us to  develop a
 vectorial Ribaucour transformation for Euclidean submanifolds,
 extending the transformation in \cite{lm} for orthogonal
 coordinate systems. It turns out that any $n$-dimensional
 submanifold with flat normal bundle of $\R^{n+m}$ can be
 obtained by applying a  vectorial Ribaucour transformation
 to an orthogonal coordinate system in an $n$-dimensional
 subspace of $\R^{n+m}$. This yields the following explicit
 local construction of {\em all\/} $n$-dimensional submanifolds
 with flat normal bundle of $\R^{n+m}$.  Notice that
carrying a principal coordinate system is not required.
\vspace{3ex}

\begin{small}
\end{small}

\newpage

\begin{theorem}\po\label{cor:ferap}
Let $\va_1,\ldots,\va_m$ be smooth real functions on
an open simply connected subset $U\subset \R^n$ satisfying
$$
[\hess \va_i,\hess\va_j]=0,\;\;\; 1\leq i,j\leq m,
$$ and let ${\cal G}\colon\,U\to
M_{n\times m}(\R)$ be defined by
${\cal G}=(\nabla\va_1,\ldots,\nabla\va_m)$.
Then for any $x\in U$ there
exists a smooth map $\Omega\colon\,V\to Gl(\R^m)$ on an
open subset $V\subset U$ containing $x$ such that
$d\Omega={\cal G}^t\,d{\cal G}\,\,\,\mbox{and}\,\,\,
\Omega+\Omega^t={\cal G}^t{\cal G}+I.$
Moreover, the map
$$
f=\left(\begin{array}{c}
id+{\cal G}\Omega^{-1}\va \\
\Omega^{-1}\va\\
\end{array}\right)
$$
with $\va^t=(\va_1,\ldots,\va_m)$ defines, at regular points, an immersion
$f\colon\,{V}\to \R^{n+m}$ with flat normal bundle.

 Conversely, any isometric immersion
$f\colon\,M^n\to \R^{n+m}$ with flat normal bundle can be locally constructed in this way.
\end{theorem}

The case of  submanifolds of the sphere can be easily derived from the
preceding result and the observation that any such submanifold arises as
the image of a unit parallel normal vector field to a submanifold with
flat normal bundle of Euclidean space (see Corollary \ref{cor:ferap2}). In
this way we recover Ferapontov's result for the holonomic case (see
Theorem \ref{thm:ferap}), thus proving his guess correct.

Theorem \ref{cor:ferap} is actually  a particular case of a more
general result. In fact, we obtain a similar local explicit
construction (see Theorem \ref{thm:parallel}) of all isometric
immersions  $\tilde{f}\colon\tilde M^{n+m}\to\R^{n+m+p}$ carrying
a parallel flat normal subbundle  of rank $m$, in particular of
all those that carry a parallel normal vector field,  starting
with an \ii $f\colon M^n\to\R^{n+p}$ and a set of Codazzi tensors
$\Phi_1,\ldots,\Phi_m$ on $M^n$ that commute one with each other
and with the second fundamental form of $f$. We refer the reader
to \cite{bco} for  results of a global nature on such isometric
immersions, with strong implications for the submanifold geometry
of orbits of orthogonal representations.

    By putting together the preceding result with
Theorem 8 of \cite{dft}, we  obtain
an explicit construction (see Theorem \ref{dupin}) in terms of
the vectorial Ribaucour transformation
of all Euclidean submanifolds that carry
a Dupin principal curvature normal vector field with integrable
conullity (see Section $7$ for the precise definitions), a concept
that has proven to be crucial
in the study of reducibility of Dupin submanifolds (see \cite{dft}).

A key feature of the Ribaucour transformation for submanifolds (in
particular,  orthogonal systems) is its permutability property. Namely,
given two Ribaucour transforms of a submanifold, there is, generically, a
fourth submanifold that is a simultaneous Ribaucour transform of the first
two, giving rise to a {\em Bianchi quadrilateral.}

     More generally, for any
integer $k\geq 2$ we define a {\em Bianchi $k$-cube\/}
as a $(k+1)$-tuple $(\mathcal{C}_0,\ldots, \mathcal{C}_k)$,
where each $\mathcal{C}_i$, $0\leq i\leq k$, is a family
of submanifolds with exactly $\binom{k}{i}$ elements,
such that every element of $\mathcal{C}_1$ is a Ribaucour
transform of the unique element of $\mathcal{C}_0$ and such
that, for every $\hat{f}\in \mathcal{C}_{s+1}$,
$1\leq s\leq k-1$, there exist unique elements
$\hat{f}_1,\ldots,\hat{f}_{s+1}\in \mathcal{C}_s$
satisfying the following conditions:
\begin{itemize}
\item[$(i)$] $\hat{f}$ is a Ribaucour transform of
$\hat{f}_1,\ldots,\hat{f}_{s+1}$.
\item[$(ii)$] For each pair of indices $1\leq i\neq j\leq s+1$
there exists a unique element
$\hat{f}_{ij}\in \mathcal{C}_{s-1}$ such that
$\{\hat{f}_{ij},\hat{f}_i,\hat{f}_j,\hat{f}\}$ is a
Bianchi quadrilateral.
\end{itemize}

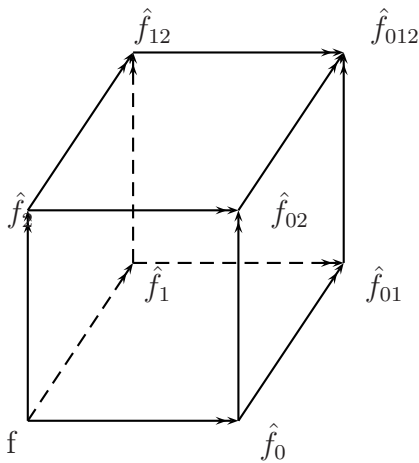
\psset{unit=7mm}
\begin{pspicture}[0](1,1)(15,10)
\psline[linestyle=dashed,arrows=->>](7,2)(9,5)
\psline[arrows=->>](7,2)(11,2)
\psline[arrows=->>](11,2)(13,5)
\psline[linestyle=dashed,arrows=->>](9,5)(13,5)
\psline[arrows=->>](7,6)(9,9)
\psline[arrows=->>](7,6)(11,6)
\psline[arrows=->>](11,6)(13,9)
\psline[arrows=->>](9,9)(13,9)
\psline[arrows=->>](7,2)(7,6)
\psline[arrows=->>](11,2)(11,6)
\psline[linestyle=dashed,arrows=->>](9,5)(9,9)
\psline[arrows=->>](13,5)(13,9)
\uput{0}[0](6.6,1.6){f}
\uput{0}[0](11.4,1.6){$\hat f_0$}
\uput{0}[0](9.2,4.6){$\hat f_1$}
\uput{0}[0](6.6,6){$\hat f_2$}
\uput{0}[0](11.6,6){$\hat f_{02}$}
\uput{0}[0](9,9.5){$\hat f_{12}$}
\uput{0}[0](13.4,4.6){$\hat f_{01}$}
\uput{0}[0](13.5,9.5){$\hat f_{012}$}
\end{pspicture}

  The following Bianchi $k$-cube theorem was proved in \cite{gt} for $k=3$ in
the context of triply orthogonal systems of Euclidean
space. A nice  proof in the setup of Lie sphere
geometry was recently given in \cite{bhj}, where also an indication
was provided of how the general
case can be settled by using results of \cite{hj} for discrete orthogonal
nets together with an induction argument.

\begin{theorem}\po\label{thm:bianchi}
Let $f\colon\,M^n\to \R^{N}$  be an isometric immersion
and let $f_1,\ldots,f_k$ be independent Ribaucour
transforms of $f$. Then, for a generic choice of  simultaneous
Ribaucour transforms $f_{ij}$ of $f_i$ and $f_j$ such that
$\{f_{ij},f_i,f_j,f\}$ is a Bianchi
quadrilateral for all pairs $\{i,j\}\subset \{1,\ldots,k\}$
with $i\neq j$, there exists a unique Bianchi $k$-cube
$(\mathcal{C}_0,\ldots, \mathcal{C}_k)$ such that
$\mathcal{C}_0=\{f\}$, $\mathcal{C}_1=\{f_1,\ldots, f_k\}$
and $\mathcal{C}_2=\{f_{ij}\}_{1\leq i\neq j\leq k}$.
\end{theorem}

We give a simple and direct proof of Theorem \ref{thm:bianchi} in
Section $5$, where the precise meanings of {\em independent\/} and
{\em generic\/} are explained. The proof relies on a general
decomposition theorem for the vectorial Ribaucour transformation
for submanifolds (Theorem \ref{viii}), according to which the
composition of two or more vectorial Ribaucour transformations
with appropriate conditions is again a vectorial Ribaucour
transformation. The latter extends a similar result of  \cite{lm}
for the case of orthogonal systems and implies, in particular, the
classical permutability of Ribaucour transformations for surfaces
and, more generally, the permutability of vectorial Ribaucour
transformations for submanifolds. \vspace{3ex}

\newpage

\noindent{\bf\large \S 1 Preliminaries.} \vspace{3ex}

 Let $M^n$ be an $n$-dimensional Riemannian manifold
and let $\xi$ be a Riemannian vector bundle over $M^n$ endowed
with a compatible connection $\nabla^\xi$. We
denote by $\Gamma(\xi)$ the space of smooth sections of $\xi$ and
by $R^\xi$ its curvature tensor. If
$\zeta=\xi^*\otimes \eta=\mbox{Hom}(\xi,\eta)$ is the
tensor product of the vector bundles $\xi^*$ and $\eta$, where $\xi^*$
stands for the dual vector bundle of $\xi$   and $\eta$ is a
Riemannian vector bundle over $M^n$, then the covariant derivative
$\nabla Z\in\Gamma(T^*M\otimes \zeta)$ of
$Z\in\Gamma(\zeta)$ is given by
$$
(\nabla^\zeta_X Z)(v)
=\nabla^\eta_X Z(v)-Z(\nabla^\xi_X v)
$$
for any $X\in \Gamma(TM)$ and  $v\in\Gamma(\xi)$. In particular,
if $\o\in \Gamma(T^*M\otimes \xi)$ is a  smooth one-form
on $M^n$ with values in $\xi$, then
$\nabla \omega\in \Gamma(T^*M\otimes T^*M\otimes \xi)$ is  given by
$$
\nabla \omega(X,Y):=(\nabla^{T^*M\otimes \xi}_X\o)(Y)=\nabla^\xi_X\omega(Y)-\omega(\nabla_X Y),
$$
where in the right hand side  $\nabla$ denotes  the Levi-Civita connection of $M^n$.
The exterior derivative $d\omega\in\Gamma(\Lambda^2T^*M\otimes \xi)$
of $\o$ is related to $\nabla\omega$ by
$$
d\omega(X,Y)=\nabla \omega(X,Y) - \nabla\omega(Y,X).
$$

 The one-form $\omega$ is {\em closed\/} if
$d\omega=0$. If $Z\in \Gamma(\xi)$, then
$\nabla Z=dZ\in\Gamma(T^*M\otimes \xi)$
is the one-form given by $\nabla
Z(X)=\nabla^\xi_X Z$. In case $\xi=M\times V$ is a trivial vector
bundle over $M^n$,  with $V$ an Euclidean vector space, that is, a
vector space endowed with an inner product, then $\Gamma(T^*M\otimes \xi)$
is identified with the space of smooth one-forms with values in $V$.
We use the same notation for the vector space $V$ and the trivial
vector bundle $\xi=M\times V$ over $M^n$.

Given $Z_1\in \Gamma(\xi^*\otimes \eta)$ and
$Z_2\in\Gamma(\eta^*\otimes\gamma)$, we define
$Z_2Z_1\in\Gamma(\xi^*\otimes\gamma)$ by
$$
Z_2Z_1(v)=Z_2(Z_1(v)),\;\; v\in \Gamma(\xi).
$$
For $Z\in \Gamma(\xi^*\otimes \eta)$, we define
$Z^t\in\Gamma(\eta^*\otimes \xi)$ by
$$
\<Z^t(u),v\>=\<u,Z(v)\>, \;\;u\in
\Gamma(\eta)\an v\in \Gamma(\xi).
$$

For later use, we summarize in the following lemma a
few elementary  properties of covariant and exterior
derivatives, which follow by straightforward
computations.

\begin{lemma}\po\label{le:prule}
The following facts
hold:\vspace{1ex}\\
$(i)$ \vspace{1ex}If $Z_1\in\Gamma(\xi^*\otimes\eta)$ and
$Z_2\in \Gamma(\eta^*\otimes\gamma)$, then $
d(Z_2Z_1)=(dZ_2)Z_1+Z_2(dZ_1)$.\\
$(ii)$  If $Z\in\Gamma(\xi^*\otimes \eta)$ then $dZ^t=(dZ)^t$.\vspace{1ex}\\
$(iii)$  If $Z\in \Gamma(\xi)$ then $ d^2Z(X,Y)=R^\xi(X,Y)Z.$ \vspace{1ex}\\
$(iv)$ If
$\zeta=\xi^*\otimes \eta$ and $Z\in \Gamma(\zeta)$ then $
(R^{\zeta}(X,Y)Z)(v)=R^\eta(X,Y)Z(v)-Z(R^\xi(X,Y)v)$.
\end{lemma}

We also need the following result.

\begin{proposition}\po\label{prop:basic}
Let $\xi, \eta$ be Riemannian vector bundles over $M^n$ and
$\omega\in \Gamma(T^*M\otimes \xi)$. Set $\zeta =\eta^*\otimes TM$
and $\gamma=\eta^*\otimes \xi$. Let $\Phi\in \Gamma(T^*M\otimes
\zeta)$ be a closed one-form such that
\be\label{eq:comut}
\nabla\o(X,\Phi_u Y)=\nabla\o(Y,\Phi_u X)
\;\;\;\mbox{for all}\;\;u\in \Gamma(\eta),
\ee
where we write
$\Phi_u X=\Phi(X)(u)$. Then the one-form $\rho=\rho(\o,\Phi)\in
\Gamma(T^*M\otimes \gamma)$ defined by
$\rho(X)(u)=\o(\Phi_uX)$
is also closed.
\end{proposition}

\proof We have
\begin{eqnarray}\label{eq:rho}
\nabla \rho(X,Y)(u)\!\!\!&=&\!\!\!\nabla^\xi_X\rho(Y)(u)
-\rho(Y)(\nabla^\eta_X u) -\rho(\nabla_X Y)(u)\nonumber\\
\!\!\!&=&\!\!\! \nabla^\xi_X \omega(\Phi_uY)
-\omega(\Phi_{\nabla^\eta_X u}Y)-\omega(\Phi_u\nabla_X Y)\nonumber\\
\!\!\!&=&\!\!\! \nabla
\omega(X,\Phi_uY)+\omega(\nabla\Phi(X,Y)u).
\,\,\,\,\,\,\qed
\vspace{1.5ex}
\end{eqnarray}

    The following consequence of Proposition \ref{prop:basic} will be used throughout the paper.

\begin{corollary}\po\label{cor:basic}
Under the assumptions of Proposition \ref{prop:basic},
assume further that $M^n$ is simply-connected and that
$\xi$ and $\eta$ are flat.
Then there exists
$\Omega(\o,\Phi)\in\Gamma(\eta^*\otimes \xi)$ such that
$$
d\Omega(\o,\Phi)(X)(u)=\o(\Phi_uX)
\,\,\,\mbox{for all}\,\,\,X\in TM
\;\;\mbox{and}\;\;u\in \eta.
$$
\end{corollary}\vspace{1ex}

\proof Since $\xi$ and $\eta$ are flat, the same holds for
$\gamma=\eta^*\otimes \xi$ by Lemma \ref{le:prule}. The manifold
$M^n$ being simply-connected,   a one-form $\rho\in
\Gamma(T^*M\otimes \gamma)$ is exact if and only if it is
closed.\vspace{3ex}\qed

\noindent{\bf\large \S 2 The Combescure transformation. } \vspace{3ex}

  In this section we introduce a vectorial version of the Combescure
transformation for submanifolds and derive a few properties of it that
will be needed later.

\begin{proposition}\po\label{prop:Vcombe0}
Let $f\colon\,M^n\to \R^N$ be an isometric
immersion of a simply connected Riemannian manifold,
let $V$ be an Euclidean vector space and let
$\Phi\in \Gamma(T^*M\otimes V^*\otimes TM)$. Then there exists
$\Fes\in\Gamma(V^*\otimes f^*T\R^N)$
such that
\be\label{eq:combe0}
d\Fes(X)(v)=f_*\Phi_vX\,\,\,\mbox{for all}\,\,\,X\in TM
\,\mbox{and}\,\,\,v\in V
\ee
if and only if $\Phi$ is closed and satisfies
\be\label{eq:comut2}
\alpha(X,\Phi_vY)=\alpha(Y,\Phi_v X)
\,\,\,\,\mbox{for all}\,\, v\in \Gamma(V),
\ee
where $\alpha\colon\,TM\times TM\to T^\perp M$
is the \sff of $f$.
\end{proposition}

\proof Applying  (\ref{eq:rho}) for $\o=f_*\in \Gamma(T^*M\otimes
f^*T\R^N)$ and $\Phi\in \Gamma(T^*M\otimes V^*\otimes TM)$
we obtain that the one-form $\rho=\rho(f_*,\Phi)\in
\Gamma(T^*M\otimes V^*\otimes f^*T\R^N)$ satisfies
$$
\nabla \rho(X,Y)u=\alpha(X,\Phi_uY)+f_*(\nabla\Phi(X,Y)u).
$$
Therefore, $\Phi$ being closed and (\ref{eq:comut2})
are both necessary and sufficient conditions for $\rho$ to be closed.
Since $V$ and $f^*T\R^N$ are flat, the result follows from  Corollary \ref{cor:basic}.
\vspace{2ex}\qed

We call ${\cal F}$  a {\it Combescure transform\/} of $f$ determined by
$\Phi$ if, in addition, \be\label{eq:sym0}
\<\Phi_vX,Y\>=\<X,\Phi_vY\>\,\,\,\,\mbox{for all}\,\,v\in \Gamma(V). \ee
Observe that ${\cal F}$ is determined up to a parallel element in
$\Gamma(V^*\oplus f^*T\R^N)$). Notice also that for each fixed vector
$v\in V$, regarded as a parallel section of the trivial vector bundle $V$,
we have that ${\cal F}(v)\in \Gamma(f^*T\R^N)$ satisfies
$d{\cal F}(v)(X)=f_*\Phi_v(X)$, and hence ${\cal F}(v)$
is a Combescure transform of $f$ in the sense of \cite{dt2}
determined by the Codazzi tensor $\Phi_v$.

\begin{proposition}\po\label{prop:Vcombe}
Let $f\colon\,M^n\to \R^N$ be an isometric immersion of a simply
connected Riemannian manifold, let $V$ be an Euclidean vector
space and let $\Phi\in \Gamma(T^*M\otimes V^*\otimes TM)$ be
closed and satisfy (\ref{eq:comut2}).  For
$\Fes\in\Gamma(V^*\otimes f^*T\R^N)$ satisfying (\ref{eq:combe0})
write \be\label{eq:OfPhi} \Fes=f_*\o^t+\beta, \ee where $\omega\in
\Gamma(T^*M\otimes V)$ and $\beta\in\Gamma(V^*\otimes T^\perp M)$.
Then \be\label{eq:alphao} \alpha(X,\o^t(v)) +(\nabla_X^{V^*\otimes
T^\perp M} \beta)v=0 \,\,\,\,\mbox{for all}\,\, v\in \Gamma(V),
\ee and $\Phi$ is given by \be\label{eq:Phi}\Phi_vX
=(\nabla^{V^*\otimes TM}_X \omega^t)v-A_{\beta(v)}X .\ee
Conversely, if $\omega\in  \Gamma(T^*M\otimes V)$ and $\beta\in
\Gamma(V^*\otimes T^\perp M)$ satisfy (\ref{eq:alphao}), then
(\ref{eq:combe0}) holds for  $\Fes=\Fes(\o,\beta)$ and
$\Phi=\Phi(\o,\beta)$ given by (\ref{eq:OfPhi}) and
(\ref{eq:Phi}), respectively. In particular, $\Phi$ is closed and
(\ref{eq:comut2}) holds. Moreover, $\Phi$ satisfies
(\ref{eq:sym0}) if and only if $\o=d\va$ for some $\va\in
\Gamma(V)$.
\end{proposition}

\proof Denote by ${\nabla}^*$ the covariant derivative of
$V^*\otimes f^*T\R^N$. Then,
\bea
d\Fes(X)(v)\!\!\!&=&\!\!\!(\nabla^*_X f_*\o^t)v
+(\nabla^*_X \beta)v=
\nabla^{f^*T\R^N}_X f_*\o^t(v) -f_*\o^t(\nabla^V_X v)+(\nabla^*_X \beta)v\\
\!\!\!&=&\!\!\! f_*\nabla_X
\o^t(v)+\alpha(X,\o^t(v))-f_*\o^t(\nabla^V_Xv)
+({\nabla}^{V^*\otimes T^\perp M}_X \beta)v-f_*A_{\beta(v)}X.
\eea
Since, on the other hand, $\Fes$ satisfies (\ref{eq:combe0}), then (\ref{eq:alphao}) and
(\ref{eq:Phi}) follow.

   Conversely, if $\o$ and $ \beta$ satisfy
(\ref{eq:alphao}), then the preceding computation  yields
(\ref{eq:combe0}) with $\Phi$  given by (\ref{eq:Phi}).
Finally, taking the inner product of (\ref{eq:Phi}) with $Y\in \Gamma(TM)$ gives
$$
\<\Phi_vX,Y\>=\<v,\nabla\o(X,Y)\>-\<A_{\beta(v)}X,Y\>,
$$
thus the symmetry of $\nabla\o$ is equivalent to
(\ref{eq:sym0}).\qed

\begin{proposition}\po\label{prop:Vcombe2} Let $f\colon\,M^n\to \R^N$ be an isometric
immersion of a simply connected Riemannian manifold. Let $V_i$, $1\leq i\leq 2$, be
Euclidean vector spaces,  and assume that
$\omega_i\in  \Gamma(T^*M\otimes V_i)$ and
$\beta_i\in \Gamma(V_i^*\otimes T^\perp M)$ satisfy
\be\label{eq:ojbj}
\alpha(X,\o_i^t(v_i))+(\nabla_X^{V_i^*\otimes
T^\perp M} \beta_i)v_i=0 \;\;\;
\mbox{for all}\;\; v_i\in\Gamma(V_i).
\ee
Set $\Fes_i=f_*\o_i^t+\beta_i$ and
$\Phi^i_{v_i}X =(\nabla^{V_i^*\otimes TM}_X
\omega_i^t)v_i-A_{\beta_i(v_i)}X$.
Then,
\be\label{eq:sym5}
\nabla \o_i(X,\Phi^j_{v_j}Y)=\nabla
\o_i(Y,\Phi^j_{v_j}X)\;\;\;\mbox{for all} \;\; v_j\in
\Gamma(V_j)
\ee
if and only if
\be\label{eq:sym6}
\<\Phi^i_{v_i}X,\Phi^j_{v_j}Y\>=
\<\Phi^i_{v_i}Y,\Phi^j_{v_j}X\>
\;\;\;\mbox{for all} \;\; v_i\in
\Gamma(V_i)\;\;\mbox{and}\;\; v_j\in \Gamma(V_j).
\ee

   When this is the case, there exists
$\Omega_{ij}
=\Omega(\o_i,\Phi^j)\in \Gamma(V_j^*\otimes V_i)$
satisfying
\be\label{eq:doij}
d\Omega_{ij}(X)(v_j)=\o_i(\Phi^j_{v_j}X)
\;\;\;\mbox{for all}\;\; v_j\in \Gamma(V_j).
\ee
In particular,
\be\label{eq:do3}
d\Omega_{ij}=\Fes_i^t\,d\Fes_j
\ee
and
\be\label{eq:OfPhi6}
\Omega_{ij}+\Omega_{ji}^t=\Fes_i^t\Fes_j
=\o_i\o_j^t+\beta_i^t\beta_j,
\ee
up to a parallel element in $\Gamma(V_j^*\otimes V_i)$.
\end{proposition}

\proof  Since $\o_j$ and $\beta_j$ satisfy (\ref{eq:ojbj}), we have
$\alpha(X,\Phi^j_{v_j}Y)=\alpha(Y,\Phi^j_{v_j}X)$ by
Proposition \ref{prop:Vcombe}. Thus, it follows from
$$
\<\Phi^i_{v_i}X,\Phi^j_{v_j}Y\>
=\<v_i,\nabla\o_i(X,\Phi^j_{v_j}Y)\>
-\<\alpha(X,\Phi^j_{v_j}Y),\beta_i(v_i)\>
$$
that conditions (\ref{eq:sym5}) and (\ref{eq:sym6}) are equivalent.
If (\ref{eq:sym5}) holds, then by Corollary \ref{cor:basic}
there exists
$\Omega_{ij}\in \Gamma(V_j^*\otimes V_i)$ satisfying (\ref{eq:doij}).  On the other hand,
$$
\Fes_i^t\,d\Fes_j(X)v_j=\Fes_i^tf_*\Phi_{v_j}^jX
=\o_i(\Phi^j_{v_j}X)
\;\;\;\mbox{for all} \;\; v_j\in \Gamma(V_j),
$$
and (\ref{eq:do3}) follows. Finally, (\ref{eq:do3})
implies that the exterior derivatives of both sides
in the first equality  of (\ref{eq:OfPhi6}) coincide.
\qed\vspace{3ex}

\noindent{\bf\large \S 3 The vectorial Ribaucour transformation.}
\vspace{3ex}

We now introduce the main concept of this paper.

\begin{definition}\po\label{ribvet}
{\em  Let $f\colon\,M^n\to \R^N$ be an isometric
immersion of a simply connected Riemannian manifold, and
let $V$  be an Euclidean  vector space. Let
$\va\in \Gamma(V)$ and
\mbox{$\beta\in \Gamma(V^*\otimes T^\perp M)$}
satisfy (\ref{eq:alphao}) with $\o=d\va$, and let
$\Omega\in\Gamma(Gl(V))$ be a solution of
the completely integrable first order system
\be\label{eq:o1}
d\Omega=\Fes^t\,d\Fes
\ee
such that
\be\label{eq:o2}
\Omega+\Omega^t=\Fes^t\Fes,
\ee
where $\Fes=f_*\o^t+\beta$. If the map
$\tilde{f}\colon\,M^n\to \R^N$ given  by
\be\label{eq:rb2}
\tilde{f}=f-\Fes\Omega^{-1}\va
\ee
is an immersion, then the isometric immersion $\tilde{f}\colon\,\tilde M^n\to \R^N$,
where $\tilde M^n$ stands for $M^n$ with the metric induced
by $\tilde f$, is called  a {\em vectorial Ribaucour transform\/}
of $f$ determined by  $(\va,\beta,\Omega)$, and it is denoted by
$\Ral_{\va,\beta,\Omega}(f)$}.
\end{definition}

\begin{remark}\po
{\em  If $\dim V=1$,  after identifying $V^*\otimes T_f^\perp M$
with $T_f^\perp M$ then $\va$ and $\beta$ become elements of $C^\infty(M)$ and
$\Gamma(T_f^\perp M)$, respectively, and
(\ref{eq:alphao}) reduces to
$$
\alpha(X,\nabla\va)+\nabla^\perp_X\beta=0.
$$
Moreover,  $\Omega=(1/2)\<\Fes,\Fes\>$ for $\Fes=f_*\nabla\va+\beta$,
and (\ref{eq:rb2}) reduces to the parameterization of a {\em scalar\/}
Ribaucour transform of $f$ obtained in Theorem $17$ of \cite{dt2}.
In this case, since $\Omega$ is determined by $\va$ and $\beta$
we write $\tilde{f}=\Ral_{\va,\beta}(f)$  instead of $\tilde{f}=\Ral_{\va,\beta,\Omega}(f)$.}
\end{remark}

Next we derive several basic properties of the vectorial Ribaucour transformation.

\begin{proposition}
\po\label{i} The bundle map $\P\in\Gamma((f^*T\R^N)^*\otimes
\tilde{f}^*T\R^N)$  given by
\be\label{eq:p}
\P=I-\Fes\Omega^{-1}\Fes^t
\ee
is a vector bundle isometry and
\be\label{eq:p1}
\tilde{f}_*=\P f_*D,\ee where
$D=I-\Phi_{\Omega^{-1}\va}\in
\Gamma(T^*M\otimes TM)$.
In particular, $\tilde{f}$ has the metric $\<\;\;,\;\>^\sim=D^*\<\;\;,\;\>.$
\end{proposition}

\proof We have
\bea
\P^t\P\!\!\!&=&\!\!\!
(I-\Fes(\Omega^{-1})^t\Fes^t)(I-\Fes\Omega^{-1}\Fes^t)\\
\!\!\!&=&\!\!\! I-\Fes\Omega^{-1}\Fes^t
-\Fes(\Omega^{-1})^t\Fes^t
+\Fes(\Omega^{-1})^t\Fes^t\Fes
\Omega^{-1}\Fes^t.
\eea
Using (\ref{eq:o2}) in the last term implies that the
three last terms cancel out. Thus ${\cal P}$ is an
isometry. Now, using
(\ref{eq:combe0})  and (\ref{eq:o1}) we obtain
\bea
\tilde{f}_*\!\!\!&=&\!\!\! f_*-d\Fes\Omega^{-1}\va+{\cal
F}\Omega^{-1}d\Omega\Omega^{-1}\va-\Fes
\Omega^{-1}\o\\\!\!\!&=&\!\!\!
f_*-f_*\Phi_{\Omega^{-1}\va}+ \Fes
\Omega^{-1}\Fes^t\,d\Fes\Omega^{-1}\va
-\Fes\Omega^{-1}\Fes^tf_*\\
\!\!\!&=&\!\!\! f_*(I-\Phi_{\Omega^{-1}\va})
-\Fes\Omega^{-1}\Fes^tf_*(I-\Phi_{\Omega^{-1}\va})
=\P f_*D. \;\;\;\vspace{1ex}
\qed
\eea

\begin{proposition}\po\label{ii}
The normal connections and second fundamental forms of $f$ and $\tilde{f}$ are related by
\be\label{eq:ncon}
\tilde{\nabla}_X^\perp¶\xi=\P\nabla_X^\perp\xi
\ee
and
\be\label{eq:sffs}
\tilde{A}_{\P\xi}=D^{-1}(A_\xi
+\Phi_{\Omega^{-1}\beta^t\xi}),
\ee
or equivalently,
\be\label{eq:sffs2}
\tilde{\alpha}(X,Y)=\P(\alpha(X,DY)
+\beta(\Omega^{-1})^t\Phi(X)^tDY).
\ee
\end{proposition}

\proof Let $\bar\nabla$ denote the connection of
$\tilde{f}^*T\R^N$. Observing that $d\Fes^t(X)$ vanishes on
$T^\perp M$, for $\<d\Fes^t(X)\xi,v\>=\<\xi,d\Fes(X)v\>
=\<\xi,f_*\Phi_v(X)\>=0$, and using (\ref{eq:combe0}),
(\ref{eq:o1}) and (\ref{eq:p1}), we get
$$
\begin{array}{l}
-\tilde{f}_*\tilde{A}_{\P\xi}X
+\tilde{\nabla}_X^\perp¶\xi
=\bar\nabla_X\P\xi=\bar\nabla_X(\xi
-\Fes\Omega^{-1}\Fes^t\xi)
\vspace{1ex}\\
=-f_*A_\xi X
+\nabla_X^\perp\xi-d\Fes(X)\Omega^{-1}\Fes^t\xi
+\Fes\Omega^{-1}d\Omega(X)\Omega^{-1}\Fes^t\xi
+\Fes\Omega^{-1}\Fes^t(f_*A_\xi X-\nabla_X^\perp\xi)
\vspace{1ex}\\
=-\P f_*A_\xi X+\P\nabla_X^\perp\xi-\P f_*\Phi(X)\Omega^{-1}\Fes^t\xi,
\end{array}
$$
which gives (\ref{eq:ncon}) and (\ref{eq:sffs}).
\qed

\begin{proposition}\po\label{iv}
The triple
$(\tilde{\va},\tilde{\beta},\tilde\Omega)=
(\Omega^{-1}\va,\P\beta(\Omega^{-1})^t,\Omega^{-1})$
satisfies the conditions of Definition \ref{ribvet} with respect
to $\tilde{f}$, and $f=\Ral_{\tilde{\va},\tilde{\beta},
\tilde\Omega}(\tilde{f})$. Moreover,
$\tilde{\Fes}=\tilde{f}_*(d\tilde{\va})^t +\tilde{\beta}$ and
$\tilde{\Phi}=\Phi(d\tilde{\va},\tilde{\beta})$ are given,
respectively, by
\be\label{eq:ftils}
\tilde{\Fes}=-\Fes\Omega^{-1}\,\an
D\tilde{\Phi}_v=-\Phi_{\Omega^{-1}v}.
\ee
\end{proposition}

\proof  Since
$\tildeø=d\tilde{\va}=-\Omega^{-1}
\o\Phi_{\Omega^{-1}\va}+\Omega^{-1}\o
=\Omega^{-1}\o D$,
we have
$$
\<\tildeø^t(v),X\>^\sim=\<v,\Omega^{-1}\o(DX)\>
=\<D\o^t(\Omega^{-1})^tv,X\>
=\<\o^t(\Omega^{-1})^tv,D^{-1}X\>^\sim,
$$
thus
\be\label{eq:ot}
D\tildeø^t=\o^t(\Omega^{-1})^t,
\ee
where we have used that $D^{-1}$ is symmetric with
respect to $\<\;,\;\>^\sim$.
We now  prove that
\be\label{eq:alpha5}
\tilde{\alpha}(X,\tildeø^t(v))+(\nabla_X^{V^*\otimes
T_{\tilde{f}}^\perp \tilde{M}} \tilde{\beta})v=0 \;\;\;\mbox{for all}\;\;\; v\in \Gamma(V).
\ee
Equations (\ref{eq:sffs2}) and (\ref{eq:ot}) yield
\be\label{eq:ta1}
\tilde{\alpha}(X,\tildeø^t(v))= \P(\alpha(X,\o^t(\Omega^{-1})^tv)
+\beta(\Omega^{-1})^t\Phi(X)^t\o^t(\Omega^{-1})^tv),
\ee
whereas (\ref{eq:ncon}) gives
\begin{eqnarray}\label{eq:ta2}
(\nabla_X^{V^*\otimes T_{\tilde{f}}^\perp \tilde{M}}
\tilde{\beta})v\!\!\!&=&\!\!\!\tilde{\nabla}_X^\perp
\tilde{\beta}(v)-\tilde{\beta}(\nabla^V_X v)=
\P({\nabla}_X^\perp\beta(\Omega^{-1})^tv
-\beta(\Omega^{-1})^t{\nabla}_X^V v)\nonumber\\
\!\!\!&=&\!\!\! \P((\nabla_X^{V^*\otimes T^\perp M} \beta)(\Omega^{-1})^tv-
\beta(\Omega^{-1})^t\Phi(X)^t\o^t(\Omega^{-1})^tv).
\end{eqnarray}
It follows from (\ref{eq:alphao}), (\ref{eq:ta1}) and
(\ref{eq:ta2}) that (\ref{eq:alpha5}) holds.

We now compute
$\tilde{\Fes}=\tilde{f}_*\tildeø^t +\tilde{\beta}$.
Using
(\ref{eq:p1}) in the first equality below, (\ref{eq:ot})
in the second and (\ref{eq:o2}) in the last one, we obtain
\be\label{eq:tildefes}
\tilde{\Fes}=\P f_*D\tildeø^t+\tilde{\beta}
=\P(f_*\o^t(\Omega^{-1})^t+\beta(\Omega^{-1})^t)
=(I-\Fes\Omega^{-1}\Fes^t)\Fes(\Omega^{-1})^t
=-\Fes\Omega^{-1}.
\ee
Then, it follows from (\ref{eq:o1}), (\ref{eq:o2}) and (\ref{eq:tildefes}) that
$$
\tilde{\Fes}^t\,d\tilde{\Fes}
=(\Omega^{-1})^t\Fes^t\,d\Fes\Omega^{-1}
-(\Omega^{-1})^t\Fes^t\Fes\Omega^{-1}d\Omega\Omega^{-1}
=d\tilde\Omega,
$$
and
$$
\tilde{\Fes}^t\tilde{\Fes}
=(\Omega^{-1})^t{\Fes}^t\Fes\Omega^{-1}
=\tilde\Omega+\tilde\Omega^t.
$$
Therefore,
$$
\Ral_{\tilde{\va},\tilde{\beta},\tilde\Omega}(\tilde{f})
=\tilde{f}-\tilde{\Fes}\tilde\Omega^{-1}\tilde{\va}
=f-\Fes\Omega^{-1}\va-(-\Fes\Omega^{-1}) \Omega\Omega^{-1}\va=f.
$$
Finally, the second formula in (\ref{eq:ftils}) follows from \bea
\tilde{f}_*\tilde{\Phi}_v(X)\!\!\!&=&\!\!\!
d\tilde{\Fes}(X)v=-d\Fes(X)\Omega^{-1}v
+\Fes\Omega^{-1}d\Omega(X)\Omega^{-1}v\\
\!\!\!&=&\!\!\! -f_*\Phi_{\Omega^{-1}v}(X)
+\Fes\Omega^{-1}\o\Phi_{\Omega^{-1}v}(X) = -\P
f_*\Phi_{\Omega^{-1}v}(X)\\\!\!\!&=&\!\!\!
-\tilde{f}_*D^{-1}\Phi_{\Omega^{-1}v}(X). \hspace{3ex}\qed\eea
\vspace{1ex}

\noindent{\bf\large \S 4 The decomposition theorem.} \vspace{3ex}

  A fundamental feature of the vectorial Ribaucour transformation for submanifolds
is the following decomposition property, first proved in
\cite{lm} in the context of orthogonal systems.

\begin{theorem}\po\label{viii}
Let
$\Ral_{\va,\beta,\Omega}(f)\colon\,\tilde{M^n}\to\R^N$
be a vectorial Ribaucour transform of an \ii
$f\colon\, M^n\to \R^N$.
For an orthogonal decomposition  \mbox{$V= V_1\oplus V_2$}  define
\be\label{vabo}
\va_j=\pi_{V_j}\circ\va,\,\,\,\,\beta_j
=\beta|_{V_j}\an
\Omega_{ij}=\pi_{V_i}\circ
\Omega|_{V_j}\in \Gamma(V_j^*\otimes V_i),
\,\,\,\;\;1\leq i,j\leq 2.
\ee
Assume that $\Omega_{jj}$ is invertible and, for $i\neq j$, define
$\Ral_{\va,\beta,\Omega}(\va_i,\beta_i,\Omega_{ii})=
(\bar\va_i,\bar\beta_i,\bar\Omega_{ii})$  by
$$
\bar\va_i=\va_i-\Omega_{ij}\Omega_{jj}^{-1}\va_j,\;\;\;
\bar\beta_i=
\P_j(\beta_i-\beta_j(\Omega^{-1}_{jj})^t\Omega_{ij}^t)\an
\bar\Omega_{ii}=
\Omega_{ii}-\Omega_{ij}\Omega_{jj}^{-1}\Omega_{ji},
$$
where $\P_j=I-\Fes_j\Omega_{jj}^{-1}\Fes_j^t$. Then the triples
$({\va}_j,{\beta}_j,{\Omega}_{jj})$ and
$(\bar\va_i,\bar\beta_i,\bar\Omega_{ii})$
satisfy the conditions of Definition \ref{ribvet} with respect to $f$ and
$f_j$, respectively, and we have
$$
\Ral_{\va,\beta,\Omega}(f)
=\Ral_{\bar\va_i,\bar\beta_i,\bar\Omega_{ii}} (\Ral_{\va_j,\beta_j,\Omega_{jj}}(f)).
$$
\end{theorem}

\proof That $({\va}_j,{\beta}_j,{\Omega}_{jj})$,
$1\leq j\leq 2$,
satisfies the conditions of Definition \ref{ribvet} with
respect to $f$ is clear. In order to prove that
$(\bar\va_i,\bar\beta_i,\bar\Omega_{ii})$ satisfies the conditions
of Definition~\ref{ribvet} with respect to
$f_j$ for $i\neq j$ we first compute $\bar\o_i=d\bar\va_i$.
We have
\bea
\bar\o_i(X)\!\!\!&=&\!\!\!\o_i(X)
+-d{\Omega}_{ij}(X)\Omega_{jj}^{-1}\va_j+
{\Omega}_{ij}\Omega_{jj}^{-1}d\Omega_{jj}(X)
\Omega_{jj}^{-1}\va_j
-{\Omega}_{ij}\Omega_{jj}^{-1}\o_j(X)\\
\!\!\!&=&\!\!\!\o_i(X)-\o_i(\Phi^j(X)
\Omega_{jj}^{-1}\va_j)
+{\Omega}_{ij}\Omega_{jj}^{-1}\o_j
(\Phi^j(X)\Omega_{jj}^{-1}\va_j)
-{\Omega}_{ij}\Omega_{jj}^{-1}\o_j(X)\\
\!\!\!&=&\!\!\!\o_i(D_jX)-{\Omega}_{ij}
\Omega_{jj}^{-1}\o_j(D_jX),
\eea
where $D_j=I-\Phi^j_{\Omega_{jj}^{-1}\va_j}$.
Then
\bea
\<\bar\o_i^t(v_i),X\>_j\!\!\!
&=&\!\!\!\<v_i,\o_i(D_jX)
-{\Omega}_{ij}\Omega_{jj}^{-1}\o_j(D_jX)\>\\
\!\!\!&=&\!\!\!\<D_j\o_i^t(v_i)
-D_j\o_j^t(\Omega_{jj}^{-1})^t{\Omega}_{ij}^t(v_i),X\>\\
\!\!\!&=&\!\!\!\<\o_i^t(v_i)
-\o_j^t(\Omega_{jj}^{-1})^t{\Omega}_{ij}^t(v_i)
,D_j^{-1}X\>_j,
\eea
where $\<\;,\;\>_j$ denotes the metric induced by $f_j$.
Using that $D_j^{-1}$ is symmetric with respect to $\<\;,\;\>_j$, we obtain that
$$
D_j\bar\o_i^t=
\o_i^t-\o_j^t(\Omega_{jj}^{-1})^t{\Omega}_{ij}^t.
$$
It follows from (\ref{eq:sffs2}) that
\be\label{eq:x}
\begin{array}{l}
\!{\alpha}_j(X,\bar\o_i^t(v_i))
\!=\P_j(\alpha(X,\o_i^t(v_i))
-\alpha(X,\o_j^t(\Omega_{jj}^{-1})^t{\Omega}_{ij}^t(v_i))
+\beta_j(\Omega_{jj}^{-1})^t\Phi^j(X)^t\o_i^t(v_i)
\vspace{1ex}\\\hspace*{13ex} +
\beta_j(\Omega_{jj}^{-1})^t
\Phi^j(X)^t\o_j^t(\Omega_{jj}^{-1})^t
{\Omega}_{ij}^t(v_i))
\end{array}
\ee
where  $\alpha_j$ is the \sff  of $f_j$.
On the other hand, we obtain from
$$
\begin{array}{l}
(\nabla_X\bar\beta_i)(v_i)
={\nabla}_X^\perp\bar\beta_i(v_i)
-\bar\beta_i(\nabla^V_Xv_i)
\vspace{1ex}\\\hspace*{7ex}
=\P_j(\nabla_X^\perp\beta_i(v_i)
-\nabla_X^\perp\beta_j(\Omega_{jj}^{-1})^t
{\Omega}_{ij}^t(v_i)-\beta_i(\nabla^V_Xv_i)
+\beta_j(\Omega_{jj}^{-1})^t{\Omega}_{ij}^t
(\nabla^V_Xv_i))
\end{array}
$$
and
$$
\begin{array}{l}
-\nabla_X^\perp\beta_j(\Omega_{jj}^{-1})^t
{\Omega}_{ij}^t(v_i)
+\beta_j(\Omega_{jj}^{-1})^t{\Omega}_{ij}^t(\nabla_Xv_i)
\vspace{1ex}\\\hspace{4ex}
=-(\nabla_X\beta_j)((\Omega_{jj}^{-1})^t
{\Omega}_{ij}^t(v_i))
-\beta_j(d(\Omega_{jj}^{-1})^t(X){\Omega}_{ij}^t(v_i))
-\beta_j(\Omega_{jj}^{-1})^t\,d{\Omega}_{ij}^t(X)(v_i)
\end{array}
$$
that
\be\label{eq:xx}\begin{array}{l}
(\nabla_X\bar\beta_i)(v_i)
=\P_j((\nabla_X\beta_i)(v_i)
-(\nabla_X\beta_j)((\Omega_{jj}^{-1})^t
{\Omega}_{ij}^t(v_i))
-\beta_j(\Omega_{jj}^{-1})^t\Phi^j(X)^t\o_i^t(v_i)
\vspace{1ex}\\\hspace{13ex}+
\beta_j(\Omega_{jj}^{-1})^t
\Phi^j(X)^t\o_j^t(\Omega_{jj}^{-1})^t
{\Omega}_{ij}^t(v_i)),
\end{array}
\ee
where we used $d{\Omega}_{ij}^t(X)=\Phi^j(X)^t\o_i^t$.
It follows from (\ref{eq:x}) and (\ref{eq:xx}) that
$$
{\alpha}_j(X,\bar\o_i^t(v_i))
+(\nabla_X\bar\beta_i)(v_i)=0.
$$
 Now we have
\be\label{eq:barfesi}
\begin{array}{l}
\bar\Fes_i={{f_j}}_{*}\bar\o_i^t+\bar\beta_i
=\P_jf_{*}D_j\bar\o_i^t+
\P_j(\beta_i-\beta_j(\Omega^{-1}_{jj})^t\Omega_{ij}^t)
\vspace{1ex}\nonumber
\\\hspace*{2.8ex}=\P_j(f_{*}\o_i^t+\beta_i
-f_{*}\o_j^t(\Omega_{jj}^{-1})^t{\Omega}_{ij}^t
-\beta_j(\Omega^{-1}_{jj})^t\Omega_{ij}^t
=\P_j(\Fes_i-\Fes_j(\Omega_{jj}^{-1})^t{\Omega}_{ij}^t)
\vspace{1ex}\nonumber\\\hspace*{2.8ex}=
\Fes_i-\Fes_j(\Omega_{jj}^{-1})^t{\Omega}_{ij}^t
-\Fes_j\Omega_{jj}^{-1}
\Fes_j^t\Fes_i+\Fes_j\Omega_{jj}^{-1}
\Fes_j^t\Fes_j(\Omega_{jj}^{-1})^t{\Omega}_{ij}^t
\vspace{1ex}\nonumber\\\hspace*{2.8ex}
=\Fes_i-\Fes_j\Omega_{jj}^{-1}\Omega_{ji},
\end{array}
\ee
where we used that $\Fes_j^t\Fes_i=\Omega_{ji}+\Omega_{ij}^t$. Then,
$$
\begin{array}{l}
\bar\Fes_i^t\,d\bar\Fes_i
=\bar\Fes_i^t(d{\Fes}_i-d{\Fes}_j
\Omega_{jj}^{-1}\Omega_{ji}
+\Fes_j\Omega_{jj}^{-1}d\Omega_{jj}\Omega_{jj}^{-1}
\Omega_{ji}
-\Fes_j\Omega_{jj}^{-1}d\Omega_{ji})
\vspace{1ex}\\\hspace*{3ex}
={\Fes}^t_id{\Fes}_i-{\Fes}^t_id{\Fes}_j
\Omega_{jj}^{-1}\Omega_{ji}
+{\Fes}^t_i{\Fes}_j\Omega_{jj}^{-1}d\Omega_{jj}
\Omega_{jj}^{-1}\Omega_{ji}
-{\Fes}^t_i{\Fes}_j\Omega_{jj}^{-1}d\Omega_{ji}
-\Omega_{ji}^t(\Omega_{jj}^{-1})^t\Fes_j^t\,d\Fes_i
\vspace{1ex}\\\hspace*{3ex}
+\Omega_{ji}^t(\Omega_{jj}^{-1})^t\Fes_j^t\,d\Fes_j
\Omega_{jj}^{-1}\Omega_{ji}
-\Omega_{ji}^t(\Omega_{jj}^{-1})^t\Fes_j^t\Fes_j
\Omega_{jj}^{-1}d\Omega_{jj}\Omega_{jj}^{-1}\Omega_{ji}
+\Omega_{ji}^t(\Omega_{jj}^{-1})^t\Fes_j^t\Fes_j
\Omega_{jj}^{-1}d\Omega_{ji}.
\end{array}
$$
Using that $d\Omega_{ji}=\Fes_j^t\,d\Fes_i$ and
$d\Omega_{jj}=\Fes_j^t\,d\Fes_j$, we obtain
$$
\bar\Fes_i^t\,d\bar\Fes_i
=d\Omega_{ii}-d\Omega_{ij}\Omega_{jj}^{-1}\Omega_{ji}
+\Omega_{ij}\Omega_{jj}^{-1}d\Omega_{jj}\Omega_{jj}^{-1}
\Omega_{ji}-\Omega_{ij}\Omega_{jj}^{-1}d\Omega_{ji}=
d\bar\Omega_{ii}.
$$
Moreover,
$$
\begin{array}{l}
\bar\Fes_i^t\bar\Fes_i=(\Fes_i^t
-\Omega_{ji}^t(\Omega_{jj}^{-1})^t\Fes_j^t)(\Fes_i
-\Fes_j\Omega_{jj}^{-1}\Omega_{ji})
\vspace{1ex}\\\hspace*{5,5ex}
=\Fes_i^t\Fes_i-\Fes_i^t\Fes_j\Omega_{jj}^{-1}\Omega_{ji}-
\Omega_{ji}^t(\Omega_{jj}^{-1})^t\Fes_j^t\Fes_i+
\Omega_{ji}^t(\Omega_{jj}^{-1})^t\Fes_j^t\Fes_j
\Omega_{jj}^{-1}\Omega_{ji}
\vspace{1ex}\\\hspace*{5,5ex}
=\Omega_{ii}-\Omega_{ij}\Omega_{jj}^{-1}
\Omega_{ji}+\Omega_{ii}^t-
\Omega_{ji}^t(\Omega_{jj}^{-1})^t\Omega_{ij}^t=
\bar\Omega_{ii}+\bar\Omega_{ii}^t,
\end{array}
$$
which completes the proof that
$(\bar\va_i,\bar\beta_i,\bar\Omega_{ii})$ satisfies the required conditions.

Now write $\Omega$ in matrix notation as
$$
\Omega=\left(\begin{array}{c}
\Omega_{11} \;\; \Omega_{12} \\
\Omega_{21} \;\; \Omega_{22}\\
\end{array}\right).
$$
Since $\Omega$ and $\Omega_{ii}$ are invertible, then
 $\bar\Omega_{ii}$ is invertible
for $1\leq i\leq 2$ and
$$
\Omega^{-1}=\left(\begin{array}{l}
\hspace*{4ex}\bar\Omega_{11}^{-1} \;\; \;\;\;\;\;\;\;\;
-\bar\Omega_{11}^{-1}\Omega_{12}\Omega_{22}^{-1}
\vspace{1ex} \\
-\bar\Omega_{22}^{-1}\Omega_{21}\Omega_{11}^{-1}
\;\;\;\;\;\;\;\;\;\;\;\; \bar\Omega_{22}^{-1}
\\
\end{array}\right).
$$
In particular,
\be\label{eq:q1}
\bar\Omega_{ii}^{-1}=\Omega_{ii}^{-1}
+\Omega_{ii}^{-1}\Omega_{ij}\bar\Omega_{jj}^{-1}
\Omega_{ji}\Omega_{ii}^{-1}
\;\an\;
\bar\Omega_{ii}^{-1}\Omega_{ij}\Omega_{jj}^{-1}=
\Omega_{ii}^{-1}\Omega_{ij}\bar\Omega_{jj}^{-1}
\ee
for $1\leq i\neq j\leq 2$. Then,
\bea
\Ral_{\va,\beta,\Omega}(f)\!\!\!&=&\!\!\!
f-\Fes\Omega^{-1}\va\\\!\!\!&=&\!\!\!f
-\Fes(\bar\Omega_{11}^{-1}(\va_1-\Omega_{12}
\Omega_{22}^{-1}\va_2)
+\bar\Omega_{22}^{-1}(\va_2-\Omega_{21}
(\Omega_{11})^{-1}\va_1))\\
\!\!\!&=&\!\!\!f-\Fes_1\bar\Omega_{11}^{-1}
(\va_1-\Omega_{12}\Omega_{22}^{-1}\va_2)
-\Fes_2\bar\Omega_{22}^{-1}(\va_2-\Omega_{21}
\Omega_{11}^{-1}\va_1).
\eea
On the other hand, by (\ref{eq:barfesi}) and
(\ref{eq:q1}) we have
\bea
\Ral_{\bar\va_i,\bar\beta_i,\bar\Omega_{ii}}({f}_j)
\!\!\!&=&\!\!\!{f}_j
-\bar\Fes_i\bar\Omega_{ii}\bar\va_i\\
\!\!\!&=&\!\!\!
f-\Fes_j\Omega_{jj}^{-1}\va_j
-(\Fes_i-\Fes_j\Omega_{jj}^{-1}\Omega_{ji})
\bar\Omega_{ii}^{-1}(\va_i
-\Omega_{ij}\Omega_{jj}^{-1}\va_j)\\
\!\!\!&=&\!\!\!
f-(\Fes_j\bar\Omega_{jj}^{-1}
-\Fes_i\bar\Omega_{ii}^{-1}\Omega_{ij}
\Omega_{jj}^{-1})\va_j
\!-\!(\Fes_i\bar\Omega_{ii}^{-1}
-\Fes_j\bar\Omega_{jj}^{-1}\Omega_{ji}
\Omega_{ii}^{-1})\va_i.
\eea
We conclude that $\Ral_{\va,\beta,\Omega}(f)
=\Ral_{\bar\va_i,\bar\beta_i,\bar\Omega_{ii}}({f}_j)$
for $1\leq i\neq j\leq 2$. \,\,\,\,\,\,\qed

\begin{remark}\po
{\em It follows from Theorem \ref{viii} that a
vectorial Ribaucour transformation whose associated data
$(\va,\beta,\Omega)$ are defined on a vector space $V$ can be
regarded as the iteration of $k=\mbox{dim} V$ scalar Ribaucour transformations.}
\end{remark}

   In applying Theorem \ref{viii}, it is often more convenient
   to use one of its two following consequences.

\begin{corollary}\po\label{x}
Let $f_i=\Ral_{\va_i,\beta_i,\Omega_{ii}}(f)\colon\, M_i^n\to\R^N$, $1\leq i\leq 2$,
be  two vectorial Ribaucour transforms of $f\colon\, M^n\to \R^N$. Assume that the tensors
${\Phi}^i=\Phi(d\va_i,{\beta}_i)$ satisfy
$$
[\Phi^i_{v_i},\Phi^j_{v_j}]=0
\;\;\;\mbox{for all}\;\;
v_i \in V_i\;\;\mbox{and}\;\;v_j\in V_j,\;\;
1\leq i\neq j\leq 2.
$$
Set $\Fes_i=f_*(d\va_i)^t+\beta_i$. Then there exists
$\Omega_{ij}\in \Gamma(V_j^*\otimes V_i)$, such that
\be\label{eq:oij1}
d\Omega_{ij}=\Fes_i^t\,d\Fes_j
\;\;\;\mbox{and}\;\;\; \Fes_i^t\Fes_j
=\Omega_{ij}+\Omega_{ji}^t,
\ee
and such that
$\va\in\Gamma(V),\beta\in\Gamma(V^*\otimes T^\perp M)$ and
$\Omega\in \Gamma(V^*\otimes V)$ defined by (\ref{vabo})
for \mbox{$V=V_1\oplus V_2$}  satisfy the conditions of
Definition \ref{ribvet} (and therefore the remaining of
the conclusions of Theorem \ref{viii} hold).
\end{corollary}

\proof The first assertion is a consequence of Proposition
\ref{prop:Vcombe2}. It is now easily seen that  $\va\in \Gamma(V),
\beta\in \Gamma(V^*\otimes T^\perp M)$, $\Omega\in
\Gamma(V^*\otimes V)$ defined by (\ref{vabo}) for $V=V_1\oplus V_2$
satisfy the  conditions of Definition \ref{ribvet} with respect
to $f$ if and only if the same holds for $(\va_i,\beta_i,\Omega_{ii})$
and, in addition, (\ref{eq:oij1}) holds.\qed

\begin{corollary}\po\label{xix} Let
$f_1=\Ral_{\va_1,\beta_1,\Omega_{11}}(f)\colon\, \tilde{M}^n\to \R^N$
be a vectorial Ribaucour transform of
$f\colon\, M^n\to \R^N$. Let $(\bar\va_2,\bar\beta_2,\bar\Omega_{22})$
satisfy the conditions of Definition \ref{ribvet} with respect to $f_1$. Assume further that
$\bar\Phi^2=\Phi(d\bar\va_2,\bar\beta_2)$
satisfies
$$
[\bar\Phi^2_{v_2},\bar\Phi^1_{v_1}]=0,
\,\,\,\mbox{for all}\,\,\,v_1\in V_1,\; v_2\in V_2,
$$
where
$D_1\bar\Phi^1_{v_1}=-\Phi^1_{\Omega_{11}^{-1}v_1}$ for
$\Phi^1=\Phi(d\va_{1},\beta_1)$. Then there exist
$\bar\Omega_{ij}\in \Gamma(V_j^*\otimes V_i)$, $i\neq j$, such that
\be\label{eq:baroij1}
d\bar\Omega_{ij}=\bar\Fes_i^t\,d\bar\Fes_j\;\;\;
\mbox{and}\;\;\;
\bar\Fes_i^t\bar\Fes_j=\bar\Omega_{ij}
+\bar\Omega_{ji}^t,
\ee
where $\bar\Fes_1=-\Fes_1\Omega_{11}^{-1}$ and
$\bar\Fes_2=(f_1)_*(d\bar\va_2)^t+\bar\beta_2$.
Now define
$$
(\va_2,\beta_2,\Omega_{22})     =(\bar\va_2-\bar\Omega_{21}\va_1,
\P_1^{-1}\bar\beta_2
-\beta_1\bar\Omega_{21}^t,
\bar\Omega_{22}
-\bar\Omega_{21}\Omega_{11}\bar\Omega_{12}),
$$
$$
\Omega_{12}=\Omega_{11}\bar\Omega_{12}
\;\;\;\mbox{and}\;\;\;
\Omega_{21}=-\bar\Omega_{21}\Omega_{11}.
$$
Then
$\va\in\Gamma(V),\beta\in\Gamma(V^*\otimes T^\perp M)$ and
$\Omega\in \Gamma(V^*\otimes V)$ defined by (\ref{vabo})
for $V=V_1\oplus V_2$  satisfy the conditions of
Definition \ref{ribvet} and
$\Ral_{\va,\beta,\Omega}(f)
=\Ral_{\bar\va_2,\bar\beta_2,\bar\Omega_{22}} (\Ral_{\va_1,\beta_1,\Omega_{11}}(f))$.
\end{corollary}

\proof By Proposition \ref{iv}, we have  $f=\Ral_{\bar\va_1,\bar\beta_1,
\bar\Omega_{11}}({f}_1)$ where
$$
(\bar\va_1,\bar\beta_1, \bar\Omega_{11})
=(\Omega_{11}^{-1}\va_1,\P_1\beta_1(\Omega_{11}^{-1})^t,
\Omega_{11}^{-1}).
$$
Moreover,
$\bar\Fes_1=(f_1)_*(d\bar\va_1)^t +\bar\beta_1$
and
$\bar\Phi_1=\Phi(d\bar\va_1,\bar\beta_1)$
are given by
$$
\bar\Fes_1=-\Fes_1\Omega_{11}^{-1}
\an D_1\bar\Phi^1_{v_1}=-\Phi^1_{\Omega_{11}^{-1}v_1}.
$$
Thus, the existence of
$\bar\Omega_{ij}\in \Gamma(V_j^*\otimes V_i)$, $i\neq j$,
satisfying condition (\ref{eq:baroij1})
follows from Proposition \ref{prop:Vcombe2} applied to
$f_1$ and the triples
$(\bar\va_i,\bar\beta_i,\bar\Omega_{ii})$,
$1\leq i\leq 2$. Now observe that
$(\va_2,\beta_2,\Omega_{22})
=\Ral_{\bar\va,\bar\beta,\bar\Omega}
(\bar\va_2,\bar\beta_2, \bar\Omega_2)$,
and hence $(\va_2,\beta_2,\Omega_{22})$ satisfies (\ref{eq:alphao})
with respect to $f_1$ and $d\Omega_{22}=\Fes_2^t\,d\Fes_2$
by Theorem \ref{viii}.

It remains to check that $d\Omega_{ij}=\Fes_i^t\,d\Fes_j$
and $\Omega_{ij}+\Omega_{ji}^t=\Fes_i^t\Fes_j$,
$\leq i\neq j$.
From the proof of
Theorem \ref{viii} (see (\ref{eq:barfesi})) we have
$\Fes_2=\bar\Fes_2-
\bar\Fes_1\bar\Omega_{11}^{-1}\bar\Omega_{12}
=\bar\Fes_2-\bar\Fes_1\Omega_{11}\bar\Omega_{12}$.
Then,
\bea
{\Fes}_2^t\,d{\Fes}_1\!\!\!&=&\!\!\!
-(\bar\Fes_2^t-\bar\Omega_{12}^t
\Omega_{11}^t\bar\Fes_1^t
\,d(\bar\Fes_1\Omega_{11})=
-\bar\Fes_2^t\,d\bar\Fes_1\Omega_{11}
-\bar\Fes_2^t\bar\Fes_1\,d\Omega_{11}+
\bar\Omega_{12}^t{\Fes}_1^t\,d{\Fes}_1\\
\!\!\!&=&\!\!\! -d\bar\Omega_{21}\Omega_{11}-(\bar\Omega_{21}
+\bar\Omega_{12}^t)\,d\Omega_{11}+
\bar\Omega_{12}^t\,d\Omega_{11}
-d(\bar\Omega_{21}\Omega_{11})\\
\!\!\!&=&\!\!\!
d\Omega_{21}.
\eea
A similar computation shows that
${\Fes}_1^t\,d{\Fes}_2=\Omega_{12}$.

Finally, we have
\bea
{\Fes}_1^t{\Fes}_2
\!\!\!&=&\!\!\!
-\Omega_{11}^t\bar\Fes_1^t(\bar\Fes_2
-\bar\Fes_1\bar\Omega_{11}^{-1}\bar\Omega_{12})
=-\Omega_{11}^t(\bar\Omega_{12}
+\bar\Omega_{21}^t-(\bar\Omega_{11}
+\bar\Omega_{11}^t)\bar\Omega_{11}^{-1}
\bar\Omega_{12})\\
\!\!\!&=&\!\!\!
-\Omega_{11}^t(\bar\Omega_{21}^t
-\bar\Omega_{11}^t\Omega_{11}\bar\Omega_{12})
=-\Omega_{11}^t\bar\Omega_{21}^t
+\Omega_{11}\bar\Omega_{12}
=\Omega_{11}\bar\Omega_{12}-
(\bar\Omega_{21}\Omega_{11})^t\\
\!\!\!&=&\!\!\!
\Omega_{12}+\Omega_{21}^t,
\eea
and similarly one checks that
${\Fes}_2^t{\Fes}_1=\Omega_{21}+\Omega_{12}^t$.
\vspace{2ex}\qed

Given four submanifolds $f_i\colon\, M_i^n\to \R^N$,
$1\leq i\leq 4$, we say that they form a {\it Bianchi
quadrilateral\/} if for each of them both the preceding
and subsequent ones (thought of as points on an oriented
circle) are Ribaucour transforms of it, and the Codazzi tensors
associated to the transformations commute.\vspace{1ex}

\noindent {\em Proof of Theorem \ref{thm:bianchi}:}
We first prove existence.  Write $f_i=\Ral_{\va_i,\beta_i}(f)$, $1\leq i\leq k$.
For each pair $\{i,j\}\subset \{1,\ldots,k\}$ with $i< j$
define $\va^{ij}\in \Gamma(\R^2)$ and $\beta^{ij}\in
\Gamma((\R^2)^*\otimes T^\perp M)$ by
$$
\va^{ij}=(\va_i,\va_j)
\;\;\;\mbox{and}\;\;\;\beta^{ij}=dx_1\otimes \beta_i+dx_2\otimes
\beta_j.
$$
By the assumption that $\{f_{ij},f_i,f_j,f\}$ is a
Bianchi quadrilateral, there
$\Omega^{ij}\in Gl(\R^2)$
with
$$
\Omega^{ij}(e_1)=\Omega_{ii}=(1/2)\<\Fes_i,\Fes_i\>\an \Omega^{ij}(e_2)
=\Omega_{jj}=(1/2)\<\Fes_j,\Fes_j\>,
$$
where
$\Fes_r=f_*\nabla\va_r+\beta_r$, $r\in\{i,j\}$,
such that
$(\va^{ij},\beta^{ij},\Omega^{ij})$ satisfies
the conditions of
Definition \ref{ribvet} with respect to $f$ and such that
$f_{ij}=\Ral_{\va^{ij},\beta^{ij},\Omega^{ij}}(f)$.
Define
$\va\in\Gamma(\R^k)$,
$\beta \in \Gamma((\R^k)^*\otimes T^\perp M)$ and
$\Omega\in \Gamma((\R^k)^*\otimes \R^k)$ by
$$
\va=(\va_1,\ldots, \va_k),\;\;\;\beta=\sum_{i=1}^k
dx_i\otimes \beta_i
$$
and
$$
\Omega= \sum_{i=1}^k \Omega_{ii} dx_i\otimes e_i +
\sum_{i<j} (\<\Omega^{ij}(e_1),e_2\> dx_i\otimes e_j
+\<\Omega^{ij}(e_2),e_1\> dx_j\otimes e_i).
$$
It is easy to check that $(\va, \beta,\Omega)$ satisfies the
conditions of Definition \ref{ribvet} with respect to $f$.

We now make precise
the ``generic" assumption on the statement of Theorem
\ref{thm:bianchi}. Namely, we require that no principal minor of
$\Omega$ vanishes, where $\Omega$ is regarded as a square
$(k\times k)$-matrix. That is, for any multi-index
$\alpha=\{i_1<\ldots <i_r\}\subset \{1,\ldots,k\}$, the sub-matrix
$\Omega_\alpha$ of $\Omega$, formed by those elements of $\Omega$
that belong to the rows and columns with indexes in $\alpha$, has
nonzero determinant. Now, for any such $\alpha$ set
$$
\va^{\alpha}=(\va_{i_1},\ldots,
\va_{i_r}),\;\;\;\beta^{\alpha}
=\sum_{j=1}^r dx_{i_j}\otimes
\beta_{i_j}\;\;\;\mbox{and}\;\;\;\Omega^\alpha
=\Omega_\alpha.
$$
We define $\mathcal{C}_r$ as the family of $\binom{k}{r}$
elements formed by the vectorial Ribaucour transforms
$\Ral_{\va^\alpha,\beta^\alpha,\Omega^\alpha}(f)$, where $\alpha$
ranges on the set of multi-indexes $\alpha=\{i_1<\ldots<i_r\}\subset \{1,\ldots,k\}$
with $r$ elements. Given
$$
\hat{f}=\Ral_{\va^\alpha,\beta^\alpha,\Omega^\alpha}(f)\in
\mathcal{C}_{s+1},\;\;\; 1\leq s\leq k-1\;\; \mbox{and} \;\;\alpha
=\{i_1<\ldots<i_{s+1}\}\subset \{1,\ldots,k\},
$$
let $\alpha_1, \ldots,\alpha_{s+1}$ be the $(s+1)$  multi-indexes  with $s$ elements
that are contained in $\alpha$. For each $j=1,\ldots, s+1$ write
$\alpha=\alpha_j\cup\{i_j\}$.  Then,
$$
\hat{f}_j:=\Ral_{\va^{\alpha_j},\beta^{\alpha_j},
\Omega^{\alpha_j}}(f)\in\mathcal{C}_s
\an
\hat{f}=\Ral_{\bar\va_{i_j},\bar\beta_{i_j}}(\hat{f}_j)
$$
by Theorem \ref{viii}. Therefore $\hat{f}$ is a Ribaucour
transform of  $\hat{f}_1,\ldots,\hat{f}_{s+1}$.
Moreover, for each pair $\{\alpha_i,\alpha_j\}$, set
$\alpha_{ij}=\alpha_i\cap \alpha_j$, and let
$\hat{f}_{ij}=\Ral_{\va^{\alpha_{ij}},
\beta^{\alpha_{ij}},\Omega^{\alpha_{ij}}}(f)$.
Then $\hat{f}_{ij}\in \mathcal{C}_{s-1}$ and
$\{\hat{f}_{ij},\hat{f}_{i},\hat{f}_{j},\hat{f}\}$ is a Bianchi quadrilateral

 Next we argue for the uniqueness. We first make precise the
meaning of $f_1,\ldots,f_k$ being {\em independent\/}  Ribaucour
transforms of $f$. Namely, if $f_i$ is determined by the pair
$(\va_i,\beta_i)$ with $\va_i\in C^\infty(M)$ and $\beta_i\in
\Gamma(T_f^\perp M)$, $1\leq i\leq k$, we require that the image of the  map
$\va=(\va_1,\ldots,\va_k)\colon\,M\to \R^k$
spans $\R^k$ and, in addition, that the linear map
\mbox{$\Fes\colon\,\R^k\to f^*T\R^N$}~given by
$\Fes=\sum_{i=1}^k dx_i\otimes \Fes_i,$
with $\Fes_i=f_*\nabla\va_i+\beta_i$, is injective.

It is easily seen that all uniqueness assertions follow from the
uniqueness for  $k=3$. For this case, the independence assumption is
equivalent to the condition that neither of $f_1$, $f_2$ or $f_3$ belong
to the {\em associated family\/} determined by the other two. Then,
uniqueness was proved in \cite{gt} by using a nice elementary argument
relying on the version of Miquel's Theorem for four circumferences.\qed
\vspace{3ex}

\noindent{\bf\large \S 5 Submanifolds carrying a
parallel flat normal subbundle.} \vspace{3ex}

In this section we give an explicit local construction
of all submanifolds of Euclidean space that carry a
parallel flat normal subbundle, from which Theorem \ref{cor:ferap} in the introduction
follows as a special case.

\begin{theorem}\po\label{thm:parallel}
Let $f\colon\,M^n\to \R^{n+p}$ be an \ii of a
simply connected Riemannian manifold and let
$\va_i\in C^\infty(M)$ and  $\beta_i\in\Gamma(T^\perp M)$,
$1\leq i\leq m$, satisfy
\be\label{eq:alpha3}
\alpha(X,\nabla\va_i)+\nabla_X^\perp\beta_i=0
\ee
and
\be\label{eq:comutfis}
[\Phi_i,\Phi_j]=0,\;\;\;1\leq i,j\leq m,
\ee
where $\Phi_i=\hess\va_i-A_{\beta_i}$.
Define ${\cal G}\colon\,M^n\to M_{(n+p)\times m}(\R)$  by
$$
{\cal G}
=(f_*\nabla\va_1+\beta_1,\ldots,f_*\nabla\va_m+\beta_m).
$$
Then there exists a smooth map
$\Omega\colon\,U\to GL(\R^m)$ on an
open subset $U\subset M^n$  such that
\be\label{eq:dv1}
d\Omega={\cal G}^t\,d{\cal G} \an
\Omega+\Omega^t={\cal G}^t{\cal G}+I.
\ee
Moreover,  the map $\tilde{f}\colon\,\tilde{M}^n\to\R^{n+p+m}$
given by
\be\label{eq:ftil}
\tilde{f}=\left(\begin{array}{c}
f-{\cal G}\Omega^{-1}\va \\
-\Omega^{-1}\va\\
\end{array}\right)
\ee
where $\va^t=(\va_1,\ldots,\va_m)$, defines, on an open subset $\tilde{M}^n\subset U$ of
regular points,  an immersion carrying a parallel flat normal subbundle of rank $m$.

 Conversely, any isometric immersion carrying a parallel flat normal subbundle
of rank~$m$ can be locally  constructed in this way.
\end{theorem}

\proof Set $V=\R^m$ and define
$\beta\in\Gamma(V^*\oplus T^\perp M)$ by $\beta(e_i)=\beta_i+e_i$,
where $\{e_i\}_{1\leq i\leq m}$
is the canonical basis of $\R^m$ regarded as the orthogonal
complement of $\R^{n+p}$ in $\R^{n+p+m}$. Then $\o=d\va$
and $\beta$ satisfy (\ref{eq:alphao})  in view of (\ref{eq:alpha3}).
Moreover, $\Fes\in \Gamma(V^*\oplus f^*T\R^{n+m+p})$
given by
$\Fes(v)=f_*\o^t(v)+\beta(v)$ satisfies
$$
\Fes=\left(\begin{array}{c}
{\cal G}\\
I_m\\
\end{array}\right),
$$
where $I_m$ denotes the $m\times m$ identity matrix. Thus
$$
\Fes^t\,d\Fes={\cal G}^t{d\cal G}\;\an\; \Fes^t\Fes
={\cal G}^t{\cal G}+I_m,
$$
and the existence of $\Omega$  satisfying (\ref{eq:dv1})
follows from  Proposition~\ref{prop:Vcombe2} by
using (\ref{eq:comutfis}). Moreover, comparing (\ref{eq:rb2}) and
(\ref{eq:ftil}) we have that
$\tilde{f}={\Ral}_{\va,\beta,\Omega}(f)$. Since $\R^m$
is a parallel flat normal vector subbundle of
$T^\perp M$ (where $f$ is
regarded as an immersion into $\R^{n+p+m}$) and
$\P\in\Gamma((f^*T\R^N)^*\otimes \tilde{f}^*T\R^N)$ given
by (\ref{eq:p}) is a parallel vector bundle isometry by virtue of
(\ref{eq:ncon}), it follows that $\P(\R^m)$ is a parallel flat
normal vector subbundle of $T_{\tilde f}M$ of rank~$m$.

 In order to prove the converse, it suffices to show that, given
 an \ii $f\colon\,M^n\to \R^{n+p+m}$ carrying a
parallel flat normal subbundle $E$ of rank $m$, there
exist locally an immersion $\tilde{f}\colon\,M^n\to\R^{n+p}\subset \R^{n+p+m}$
and $\tilde{\va}\in \Gamma(V:=\R^m)$,
$\tilde{\beta}\in\Gamma(V^*\oplus T_{\tilde{f}}^\perp M)$ and
$\tilde\Omega\in \Gamma(Gl(V))$ satisfying the conditions of
Definition \ref{ribvet} such that
$(\tilde{\beta}(e_i))_{\R^m}=e_i$, $1\leq i\leq m$, and
that
$f={\Ral}_{\tilde{\va},\tilde{\beta},\tilde\Omega}
(\tilde f)$.

 Let $\xi_1,\ldots,\xi_m$ be an orthonormal parallel
frame of $E$. Let $V=\R^m$ be identified with a subspace
of $\R^N$ and let $e_1,\ldots,e_m$ be the canonical basis
of $\R^m$. Define $\va\in\Gamma(V)$ and
$\beta\in \Gamma(V^*\otimes T^\perp M)$ by
$$
\va=-\sum_{i=1}^m\<f,e_i\>e_i\;\;\;\mbox{and}
\;\;\;\beta(v)=\sum_{i=1}^mx_i(\xi_i-e_i^\perp),
$$
for $v=(x_1,\ldots,x_m)$, where $e_i^\perp$ denotes the normal
vector field obtained by orthogonally projecting $e_i$ pointwise
onto $T^\perp M$. Then
\be\label{eq:ob}
\o^t(v)=-\sum_{i=1}^mx_if_*^te_i\;\;\;\mbox{and}\;\;\;
\Fes(v)=\sum_{i=1}^mx_i(\xi_i-e_i).
\ee
Therefore,
\bea
\alpha(X,\o^t(v))
+(\nabla_X^{V^*\otimes T^\perp M}\beta)v
\!\!\!&=&\!\!\!-\sum_{i=1}^mx_i\alpha(X,f_*^te_i)+
\sum_{i=1}^mx_i\nabla_X^\perp(\xi_i- e_i^\perp)\\
\!\!\!&=&\!\!\!-\sum_{i=1}^mx_i(\tilde\nabla_X
(f_*^te_i+e_i^\perp))^\perp=
\sum_{i=1}^mx_i(\tilde{\nabla}_X e_i)^\perp=0,
\eea
where $\tilde\nabla$ denotes the Euclidean connection, and hence (\ref{eq:alphao}) is satisfied.

It also follows from (\ref{eq:ob}) that
$$
\<\Fes^t\xi_j,v\>=\<\xi_j,\Fes(v)\>=x_j
-\sum_{i=1}^mx_i\<\xi_j,e_i\>,
$$
thus
$\Fes^t\xi_j=e_j-\sum_{i=1}^m\<\xi_j,e_i\>e_i$.
Similarly,
$\Fes^te_j=-e_j+\sum_{i=1}^m\<\xi_i,e_j\>e_i$.
We obtain,
$$
\Fes^t\Fes(v)=\sum_{j=1}^mx_j\Fes^t(\xi_j-e_j)
=\sum_{j=1}^m2x_je_j-\sum_{i,j=1}^mx_j(\<\xi_j,e_i\>
+\<\xi_i,e_j\>)e_i.
$$
In matrix notation, this reads as
$$
\Fes^t\Fes=2I-(\<\xi_j,e_i\>)-(\<\xi_i,e_j\>).
$$
Therefore $\Omega=I-(\<\xi_i,e_j\>)$ satisfies (\ref{eq:o1}) and (\ref{eq:o2}).
Moreover, since
$$
\Omega e_j=e_j-\sum_{i=1}^m\<\xi_j,e_i\>e_i=\Fes^t\xi_j,
$$
we have
$$
\P\xi_j=\xi_j-\Fes\Omega^{-1}\Fes^t\xi_j=\xi_j-\Fes
e_j=\xi_j-(\xi_j-e_j)=e_j.
$$
Therefore $\tilde{f}={\Ral}_{\va,\beta,\Omega}(f)$ is such that
$\tilde{f}(M^n)$ is contained in an affine subspace orthogonal to
$\R^m$. Since $f={\Ral}_{\tilde{\va},\tilde{\beta},
\tilde\Omega}(\tilde f)$ with the triple
$(\tilde{\va},\tilde{\beta},\tilde\Omega)$ given by Proposition
\ref{iv}, in order to complete the proof of the theorem it remains
to show that $(\tilde{\beta}(e_i))_{\R^m}=e_i$. But this follows
from
$$
\<\tilde{\beta}e_i,e_j\>
=\<\P\beta(\Omega^{-1})^te_i,e_j\>=
\<e_i,\Omega^{-1}\Fes^t\P^t e_j\>
=\<e_i,\Omega^{-1}\Fes^t\xi_j\>=\<e_i,e_j\>.
\,\,\,\,\,\,\qed\vspace{1ex}
$$

 The case of submanifolds with flat normal bundle of the sphere
 now follows easily from Theorem  \ref{cor:ferap}.

\begin{corollary}\po\label{cor:ferap2} Let $U\subset \R^n$,
$\{\va_i\}_{1\leq  i\leq m}$, ${\cal G}\colon\,U\to M_{n\times m}(\R)$,
and  $\Omega\colon\,V\subset U\to Gl(\R^m)$ be as in
Theorem \ref{cor:ferap}. Then the
$M_{(n+m)\times m}(\R)$-valued map
$$
W=\left(\begin{array}{c}
{\cal G}\Omega^{-1} \\
I-\Omega^{-1}\\
\end{array}\right)
$$
satisfies $W^tW=I$ and any of its columns defines, at regular points,
the position vector of an immersion with flat normal bundle into $\Sf^{n+m-1}$.

Conversely, any isometric immersion with flat normal bundle $f\colon\,M^n\to\Sf^{n+m-1}$
can be locally constructed in this way.
\end{corollary}
\proof Set $V=\R^m$ and define
$\beta\in\Gamma(V^*\oplus T^\perp M)$
by $\beta(e_i)=e_i$,  where $\{e_i\}_{1\leq i\leq m}$
is the canonical basis of $\R^m$ regarded as the
orthogonal complement of
$\R^{n+p}$ in $\R^{n+p+m}$. Then $\o=d\va$ and $\beta$ trivially satisfy (\ref{eq:alphao}).
Moreover, if
$\Fes\in \Gamma(V^*\oplus f^*T\R^{n+m+p})$
is given by $\Fes(v)=f_*\o^t(v)+\beta(v)$ then
$$
\Fes=\left(\begin{array}{c}
{\cal G}\\
I_m\\
\end{array}\right),
$$
where $I_m$ denotes the $m\times m$ identity matrix.
Let $\tilde{f}={\Ral}_{\va,\beta,\Omega}(id)$,
where $id$ is the
inclusion of $U$ into $\R^n$, Then the isometry $\P$ as in
(\ref{eq:p}) is given by
$$
\left(\begin{array}{c}
I-{\cal G}\Omega^{-1}{\cal G}^t \,\,\,\,\,\,\,
{\cal G}\Omega^{-1}   \\
\,\,\,\,\,\Omega^{-1}{\cal G}^t
\,\,\,\,\,\,\,\,\,\,\,\,\,I-\Omega^{-1}\\
\end{array}\right).
$$
Therefore $W^tW=I$ and the $(n+p+j)^{th}$-column of $W$ is $\P e_j$, $1\leq j\leq m$.
Therefore it is a unit parallel normal
vector field to $f$, and hence defines, at regular points,
the position vector of a submanifold with flat normal
bundle of $\Sf^{n+m-1}$.

  The converse follows from the converse in Theorem \ref{cor:ferap} and the fact that any
isometric immersion $f\colon\,M^n\to \Sf^{n+m-1}$ arises
as a parallel unit normal vector field of an isometric immersion $F\colon\,M^n\to \R^{n+m}$,
for instance, $F=i\circ f$, where $i$ is the canonical inclusion of $\Sf^{n+m-1}$ into $\R^{n+m}$.
\vspace{2ex}\qed

 We now give a precise statement of Ferapontov's theorem referred to in the introduction
for the case of  holonomic submanifolds
of the sphere, and show how it can be derived from Corollary \ref{cor:ferap2}.

\begin{theorem}\po\label{thm:ferap}
On an open simply connected subset $U\subset \R^n$
let $\{\beta_{ij}\}_{1\leq i\neq j\leq n}$ be smooth real functions satisfying the
completely integrable system of PDE's
\be\label{eq:I}
\left\{\begin{array}{l}
{\displaystyle\frac{\d \beta_{ij}}{\d u_k}
= \beta_{ik}\beta_{kj}},\;\;\; 1\leq i\neq j\neq k\neq i\leq n,
\vspace{1.5ex}\\
{\displaystyle \frac{\d \beta_{ij}}{\d u_i} + \frac{\d
\beta_{ji}}{\d u_j}\! +
\!\sum_{k}\beta_{ki}\beta_{kj}=0},\;\;i\neq j,
\vspace{1.5ex}\\
\end{array}\right.
\ee let $H^\alpha=(H_1^\alpha,\ldots,H_n^\alpha)$,
$1\leq\alpha\leq m$, be arbitrary solutions of the linear system
of PDE's \be\label{eq:hi} \frac{\d H_j}{\d u_i} =
\beta_{ij}H_i,\;\;\;1\leq i\neq j\leq n, \ee and let
$X_i\colon\,U\to\R^n$, $1\leq i\leq n$, satisfy \be\label{eq:xi}
\frac{\d X_{i}}{\d u_j} = \beta_{ij}X_j,\,\,\, i\neq j, \;\;\;\;
\frac{\d X_i}{\d u_i} =-\sum_{k\neq i}\beta_{ki}X_k \ee and
$X^tX=I$ at some point of $U$, where $X=(X_1,\ldots,X_n)$, the
integrability conditions of (\ref{eq:hi}) and (\ref{eq:xi}) being
satisfied by virtue of (\ref{eq:I}). Then there exist vector
functions $s^\alpha= (s_1^\alpha,\ldots,s_n^\alpha)\colon\,U\to
\R^n$, $1\leq\alpha\leq m$, such that
$ds_i^\alpha=\sum_{k=1}^nX_{ik}H_k^\alpha du^k$, and a
map~\mbox{$\Omega\colon\,U\to M_{n\times m}(\R)$} such that
$$
d\Omega={\cal G}^t\,d{\cal
G}\;\;\;\mbox{and}\;\;\;\Omega+\Omega^t={\cal G}^t{\cal G}+I_m,
$$
where ${\cal G}=(s^1,...,s^m)$. Moreover, the $(m+n)\times
m$-matrix
$$
W=\left(\begin{array}{ccc}
\,\,\,{\cal G}\Omega^{-1}  \\
\,\cdots\\
\,\,\,I_m-\Omega^{-1}  \\
\end{array}\right)
$$
satisfies $W^tW=I_m$ and any of its columns defines, at
regular points,  the position vector of an $n$-dimensional
submanifold $M^n\subset \Sf^{n+m-1}\subset \R^{n+m}$ with
flat normal bundle such that $u_1,\ldots,u_n$ are principal coordinates of $M^n$.

Conversely, any $n$-dimensional submanifold with flat normal
bundle of $\Sf^{n+m-1}$ carrying a holonomic net of
curvature lines can be locally constructed in this way.
\end{theorem}

\proof  It is easily checked using (\ref{eq:xi}), and the fact
that $X^tX=I$ at some point of $U$, that $X^tX=I$ everywhere on
$U$, whence $X_1,\ldots, X_n$ determine an orthonormal frame on
$U$. Define $\Phi^\alpha\in \Gamma(T^*U\otimes TU)$ by
$$
\Phi^\alpha X_i=H^\alpha_iX_i,\,\,\,1\leq i\leq n,\,\,\,1
\leq \alpha\leq m.
$$
Then $\Phi^\alpha$ is a symmetric tensor and, by (\ref{eq:hi}) and
(\ref{eq:xi}),
$$
\frac{\partial}{\partial u_j}(\Phi^\alpha X_i)
=\frac{\partial H^\alpha_i}{\partial u_j}X_i
+H^\alpha_i\frac{\partial X_i}{\partial u_j}
=\beta_{ji}H_j^\alpha X_i+H^\alpha_i\beta_{ij}X_j
=\frac{\partial}{\partial u_i}(\Phi^\alpha X_j),
$$
hence $\Phi^\alpha$ is a Codazzi tensor on $U$.
Thus $\Phi^\alpha$ is closed as a one-form in $U$ with values
in $TU$. Since $U$ is flat, there exists $Z^\alpha\in \Gamma(TU)$
such that $\Phi^\alpha=dZ^\alpha$. Moreover,
the symmetry of $\Phi^\alpha$ implies that
$Z^\alpha=\grad\, \va^\alpha$ for some
$\va^\alpha\in C^\infty(U)$, and hence $\Phi^\alpha=\mbox{\em\hess}\,\va^\alpha$
(cf.\ \cite{f}). Since
$\{X_i\}_{1\leq i\leq n}$ is a common diagonalizing
basis of $\mbox{\em\hess} \,\va^\alpha$,
$1\leq \alpha\leq m$,  it follows that
$[\mbox{\em\hess} \,\va_\alpha,\mbox{\em\hess}\,
\va_\beta]=0$, $1\leq \alpha,\beta\leq m$. Setting
$s^\alpha=\grad\,\varphi^\alpha$,
the remaining of the proof follows from Corollary
\ref{cor:ferap2}.\qed

\begin{remark}\po{\em Equations (\ref{eq:I}) (Lam\'e equations)
and (\ref{eq:hi}) are well-known in the theory of $n$-orthogonal systems
(cf.\ \cite{da}) where the functions $H_i$ and $\beta_{ij}$ are usually
called the Lam\'e and rotation coefficients, respectively. }\vspace{1ex}
\end{remark}

\newpage

\noindent{\bf\large \S 6 Submanifolds carrying a Dupin principal normal. }
\vspace{3ex}

   A smooth normal vector field $\eta$ of an isometric immersion $f\colon\,M^n\to \R^N$ is called a
{\em principal normal\/} with multiplicity $m\geq 1$ if
the tangent subspaces
$$
\E_\eta=\ker\, (\alpha-\<\,\,\,,\,\,\>\,\eta)
$$
have constant dimension $m\geq 1$.
If $\eta$ is parallel in the normal connection
along the  {\em nullity\/} distribution $\E_\eta$ then $\eta$
is said to be a {\em Dupin\/} principal normal.
This condition is automatic if $m\geq 2$. If $\eta$ is nowhere
vanishing, it is well-known that $\E_\eta$ is an involutive
distribution whose leaves are round $m$-dimensional spheres in $\R^N$. When $\eta$
vanishes identically, the distribution $\E_\eta=\E_0$ is known
as the relative nullity distribution, in which case the
leaves are open subsets of affine subspaces of $\R^N$.

    Let $h\colon\,L^{n-m}\to\R^{N}$ be an isometric
immersion carrying a parallel flat normal subbundle $\N$
of rank $m$, and let $\va\in C^\infty(L^{n-m})$ and
$\beta\in\Gamma(T^\perp_hL)$ satisfy
$$
\alpha(X,\nabla\va)+\nabla^\perp_X\beta=0.
$$
Assume further that the tangent subspaces
$$
E(x)=\{Z\in T_xL:(\alpha(Z,X))_{\N^\perp}
=\va^{-1}\beta_{\N^\perp}\<Z,X\>\,\,
\mbox{for all}\,\,X\in T_xL\}
$$
are everywhere trivial. Define
$\Ral^\N_{\va,\beta}(h)\colon\,\N\to\R^N$
by
$$
\Ral^\N_{\va,\beta}(h)(t)
=\Ral_{\va,\beta+t^{'}}(h)(x),
$$
where $x=\pi(t)$ and $t^{'}$ is the parallel section in
$\N$ such that $t^{'}(x)=t$. It was shown in  \cite{dft}
that $\Ral^\N_{\va,\beta}(h)$ defines, at regular points,
an immersion carrying a Dupin principal normal with
integrable {\em conullity} distribution $\E_\eta^\perp$
and that, conversely, any such immersion
can be~locally constructed in this way.

    Using the results of the previous sections we now give
an explicit description of all isometric immersions carrying a Dupin
principal normal of multiplicity $m$ and integrable conullity  in terms of
the vectorial Ribaucour transformation, starting with an isometric
immersion $g\colon\,L^{n-m}\to\R^{N}$ such that $g(L^{n-m})$ lies in an
$(N-m)$-dimensional subspace $\R^{N-m}\subset \R^N$.

Namely, let $\Ral_{\va, \beta,\Omega}(g)$
be a vectorial Ribaucour transform of an isometric
immersion $g\colon\,L^{n-m}\to\R^{N-m}\subset \R^N$ determined by
$(\va,\beta,\Omega)$ as in Definition \ref{ribvet}. For an orthogonal
decomposition $\R^{m+1}=\R\oplus \R^m$, with $\R=\spa\{e_0\}$, set
$$
(\va_0,\beta_0)=(\<\va,e_0\>,\beta(e_0)),\;\;\;
(\va_1,\beta_1)=(\pi_{\R^m}\circ \va,\beta|_{\R^m}),
$$
and
$$
(\bar\va_0,\hat{\beta}_0)=  (\va_0-\Omega_{01}\Omega_{11}^{-1}\va_1,
\beta_0-\beta_1(\Omega_{11}^{-1})^t\Omega_{01}^te_0)
$$
where
$$
\Omega_{11}=\pi_{\R^m}\circ \Omega|_{\R^m} \an
\Omega_{01}=\pi_{\R}\circ\Omega|_{\R^m}.$$
Assume that the bilinear maps $\gamma\colon\,T_xL\times T_xL\to T_xL^\perp$ given by
$$
\gamma(Z,X)=(\alpha_g(Z,D_1X)+
\beta_1(\Omega_{11}^{-1})^t\Phi_1(Z)D_1(X)
-\bar\va_0\<Z,X\>\hat{\beta}_0)_{\R^{N-m}}
$$
have everywhere trivial kernel. Let the subspace $\R^m=e_0^\perp$
be identified with the orthogonal complement of $\R^{N-m}$ in $\R^N$ and choose an
orthonormal basis $\{e_0,\ldots, e_m\}$ of
$V:=\R^{m+1}$. Finally, for $t=\sum_{i=1}^mt_ie_i\in \R^m$ define
$$
\beta_t
=e_0^*\otimes (\beta_0+t)+\sum_{i=1}^m
e_i^*\otimes\beta(e_i) \an\Omega_t
=\Omega+
(\<\beta_0,t\>+(1/2)|t|^2)e_0^*\otimes e_0.
$$
\begin{theorem}\po\label{dupin}  The triple
$(\va,\beta_t,\Omega_t)$ satisfies the conditions of Definition
\ref{ribvet} with respect to $g$ for each $t\in \R^m$
and the map
$G\colon\, L^{n-m}\times \R^m\to \R^N$ given by
$$
G(x,t)=\Ral_{\va,\beta_t,\Omega_t}(g)(x)
$$
parameterizes, at regular points, an $n$-dimensional
submanifold carrying a Dupin principal normal of multiplicity
$m$ with integrable conullity.

 Conversely, any isometric immersion carrying a Dupin
principal normal of multiplicity $m$ with integrable conullity
can be locally constructed in this way.
\end{theorem}

\proof The
first assertion is easily checked. By Theorem \ref{viii}
we have
$$
\Ral_{\va,\beta_t,\Omega_t}(g)=
\Ral_{\bar\va_0,\bar\beta^0_t}
(\Ral_{\va_1,\beta_1,\Omega_{11}}(g)),
$$
where
$$
\bar\beta^0_t=\P_1(\beta_0
+t-\beta_1(\Omega_{11}^{-1})^t\Omega_{01}^te_0)
=\bar\beta^0
+\sum_{i=1}^mt_i\eta_i,
$$
with $\bar\beta^0=\P_1\hat{\beta}_0$
and  $\eta_i=\P_1e_i$. Then $\N=\P_1\R^m$ is
a parallel flat
normal subbundle $\N$ of rank $m$ of
$h=\Ral_{\va_1,\beta_1,\Omega_{11}}(g)$ and
$$
G(x,t)=\Ral_{\bar\va_0,\bar\beta_0+t}(h)(x)=
\Ral^\N_{\bar\va_0,\bar\beta_0}(h)(t),
$$
where $t=\sum_{i=1}^mt_i\eta_i$. Moreover, the assumption
on the bilinear map $\gamma$ is easily seen
to be equivalent to the subspaces
$$
E(x)=\{Z\in T_xL :
(\alpha_h(Z,X))_{\N^\perp}
=-\bar\va_0^{-1}\bar\beta_{\N^\perp}\<Z,X\>
\,\,\mbox{for all}\,\,X\in T_xL\}
$$
being everywhere trivial. By the result of
\cite{dft} discussed
before the statement of Theorem~\ref{dupin}, it
follows that $G$ parameterizes, at regular points, an
$n$-dimensional submanifold carrying a Dupin principal normal
of multiplicity $m$ with integrable conullity.

Conversely, given a submanifold $f\colon\,M^m\to \R^N$
that carries a Dupin principal normal of multiplicity
$m$ with integrable conullity, by the aforementioned
result of \cite{dft} there exist an isometric immersion $h\colon\,L^{n-m}\to\R^{N}$
carrying a parallel flat normal subbundle $\N$ of rank
$m$, $\va\in C^\infty(L^{n-m})$ and
$\beta\in \Gamma(T^\perp_hL)$
satisfying $\alpha(X,\nabla\va)+\nabla^\perp_X\beta=0$,  with
$$
E(x)=\{Z\in T_xL :(\alpha(Z,X))_{\N^\perp}
=-\va^{-1}\beta_{\N^\perp}\<Z,X\>\,\,
\mbox{for all}\,\,X\in T_xL\}
$$
everywhere trivial, such that $f$  is
parameterized  by the map
$\Ral^\N_{\va,\beta}(h)\colon\,L^{n-m}\times \R^m\to\R^N$ given by
$$
\Ral^\N_{\va,\beta}(h)(x,t)
=\Ral_{\va,\beta+t^{'}}(h)(x),
$$
where $t^{'}=\sum_{i=1}^mt_i\eta_i$ for some
orthonormal parallel
frame $\eta_1,\ldots, \eta_m$ of $\N$.

As in the proof of Theorem \ref{thm:parallel} there is an \ii
$g\colon\,L^{n-m}\to\R^{N-m}\subset \R^N$ such that
$h=\Ral_{\va_1,\beta_1,\Omega_{11}}(g)$, and hence
$$
\Ral^\N_{\va,\beta}(h)(x,t)
=\Ral_{\va,\beta+t^{'}}
(\Ral_{\va_1,\beta_1,\Omega_{11}}(g))(x).
$$

    In order to apply Corollary
\ref{xix}, we must verify that the tensor
$\Phi=\mbox{Hess}\,\va-A_\beta$
associated to $(\va,\beta)$ commutes with $\tilde{\Phi}_{v_1}$ for every
$v_1\in V_1$, where $D_1\bar\Phi^1_{v_1}=-\Phi^1_{\Omega_{11}^{-1}v_1}$ for
$\Phi^1=\Phi(d\va_{1},\beta_1)$. Since  $\Phi$ commutes
with the shape operator $\tilde{A}_{\tilde{\xi}}$ of $h$ with respect to any normal vector field
$\tilde{\xi}\in \Gamma(T_{h}^\perp L)$,
it commutes in particular  with $\tilde{A}_{\P_1e_i}$.
But by (\ref{eq:sffs}) we have
$$
\tilde{A}_{\P_1e_i}
=D_1^{-1}\Phi^1_{\Omega_{11}^{-1}\beta_1^te_i}
=D_1^{-1}\Phi^1_{\Omega_{11}^{-1}e_i}
=-\bar\Phi^1_{e_i},
$$
and we are done.
It follows from Corollary \ref{xix} that there exist
$(\va,\beta,\Omega)$ satisfying the~conditions of
Definition \ref{ribvet} with respect to $g$ and an orthogonal decomposition
$\R^{m+1}=\R\oplus \R^m$ such that
$$
(\va_1,\beta_1,\Omega_{11})=(\pi_{\R^m}\circ \va,
\beta|_{\R^m},\pi_{\R^m}\circ\Omega|_{\R^m})
$$
and, setting $\va_0=\pi_{\R}\circ \va$,
$\beta_1=\beta|_{\R}$, then
$(\va,\beta)=\Ral_{\va,\beta,\Omega}(\va_0,\beta_0)$ and
$$
\Ral_{\va,\beta}(\Ral_{\va_1,\beta_1,\Omega_{11}}(g))
=\Ral_{\va,\beta,\Omega}(g).
$$
Defining $\beta_t$ and $\Omega_t$ as in the statement, we have
$\Ral_{\va,\beta+t^{'}}
(\Ral_{\va_1,\beta_1,\Omega_{11}}(g))
=\Ral_{\va,\beta_t,\Omega_t}(g)$.\qed

{\renewcommand{\baselinestretch}{1}
\hspace*{-20ex}\begin{tabbing}
\indent \= IMPA  \hspace{30ex} Universidade Federal de S\~ao Carlos \\
\>  Estrada Dona Castorina, 110 \hspace{7ex}
Via Washington Luiz km 235 \\
\> 22460-320 --- Rio de Janeiro
\hspace{8ex} 13565-905 --- S\~ao Carlos  \\
\> Brazil\hspace{31ex} Brazil\\
\> marcos@impa.br,\,\, luis@impa.br  \hspace{5ex}
tojeiro@dm.ufscar.br
\end{tabbing}}


\begin{thebibliography}{l}

\bibitem{bco} J.Berndt, S. Console and C. Olmos, Submanifolds and Holonomy, CRC/Chapman and Hall
Research Notes Series in Mathematics {\bf 434} (2003), Boca
Ratton.


\bibitem{bhj} F.E. Burstall and U. Hertrich-Jeromin,
{\it The Ribaucour transformation in Lie
sphere geometry \/},
Preprint arXiv:math.DG/$0407244$ v$1$, $2004$.

\bibitem{da} {\it G. Darboux},  Le\c cons sur les
systèmes orthogonaux et les coordonnées curvilignes, Gauthier-Villars, Paris, 1910.

\bibitem{dt1}  M. Dajczer and R. Tojeiro,
{\it An extension of the classical Ribaucour transformation.} Proc. London Math. Soc. {\bf 85} (2002),
211--232.

\bibitem{dt2} M. Dajczer and R. Tojeiro,
{\it Commuting Codazzi tensors and the Ribaucour
transformations for submanifolds.}  Result. Math. {\bf 44} (2003), 258--278.

\bibitem{dft} M. Dajczer, L. Florit and R. Tojeiro, {\it Reducibility of Dupin submanifolds.}
To appear in Illinois J. Math.

\bibitem{f} D. Ferus,  {\it A remark on Codazzi
tensors in constant curvature spaces.} Lect.
Notes in Math. 838, Berlin 1981, 247.

\bibitem{fe} E. Ferapontov,
{\it Surfaces with flat normal bundle: an explicit construction.}
Diff. Geom. Appl. {\bf 14} (2001), 15--37.

\bibitem{gt} E. Ghanza and S. Tsarev,
{\it An algebraic superposition formula and the completeness of
Bäcklund transformations of $(2+1)$-dimensional integrable systems.}
Uspekhi Mat. Nauk {\bf 51} (1996),  197--198; translation in Russian
Math. Surveys {\bf 51} (1996), 1200--1202.


\bibitem{hj} U. Hertrich-Jeromin, {\it Introduction to M\"obius differential
geometry,\/} London Mathematical Society Lecture Note Series, vol. {\bf 300},
Cambridge University Press, Cambridge, $2003$.

\bibitem{lm} Q. P. Liu and M. Manas,
{\it Vectorial Ribaucour transformation for the Lamé equations.} J. Phys. A. {\bf 31}
(1998), 193--200.



\end{thebibliography}
\end{document}